\documentclass[a4paper]{article}
\usepackage{geometry}
\usepackage[colorlinks=true,citecolor=blue]{hyperref}
\usepackage[dvips]{epsfig}
\usepackage{color}
\usepackage{bm}
\usepackage{amssymb}
\usepackage{amsmath}
\usepackage{amsfonts}
\usepackage{subfigure}
\usepackage{xspace}
\usepackage{upgreek}
\usepackage{graphicx}
\usepackage{epstopdf}
\usepackage{float}
\usepackage[table]{xcolor}
\begin{document}

\title{Radial basis function ENO and WENO finite difference methods based on the
optimization of shape parameters}


\author{Jingyang Guo \and Jae-Hun Jung
}

\author{Jingyang Guo\thanks{Department of Mathematics, University at Buffalo, SUNY, Buffalo, NY 14260-2900, USA. ({\tt jguo4@buffalo.edu}).} \and 
            Jae-Hun Jung\thanks{Corresponding author. Department of Mathematics, University at Buffalo, SUNY, Buffalo, NY 14260-2900, USA. ({\tt jaehun@buffalo.edu}).} }


\maketitle
\noindent

\section*{Abstract}
We present adaptive finite difference ENO/WENO methods by adopting infinitely smooth radial basis functions (RBFs). This is a direct extension of the non-polynomial finite volume ENO/WENO method proposed by authors in \cite{GuoJung} to the finite difference ENO/WENO method based on the original smoothness indicator scheme developed by Jiang and Shu \cite{WENO}. The RBF-ENO/WENO finite difference method slightly perturbs the reconstruction coefficients with RBFs as the reconstruction basis and enhances accuracy in the smooth region by locally optimizing the shape parameters. The RBF-ENO/WENO finite difference methods provide more accurate reconstruction than the regular ENO/WENO reconstruction and provide sharper solution profiles near the jump discontinuity. Furthermore the RBF-ENO/WENO methods are easy to implement in the existing regular ENO/WENO code. The numerical results in 1D and 2D presented in this work show that the proposed RBF-ENO/WENO finite difference method better performs than the regular ENO/WENO method.

\textbf{keywords} Essentially non-oscillatory method, Weighted essentially non-oscillatory method, Radial basis function interpolation, Finite difference method.




\section{Introduction}
We are interested in solving hyperbolic problems that could contain discontinuous solutions in time. High order numerical approximations of such discontinuous solutions may suffer from the Gibbs phenomenon and the approximations are highly oscillatory, e.g. spectral approximations of discontinuous functions \cite{Hesthaven}. High order essentially non-oscillatory (ENO) and weighted essentially non-oscillatory (WENO) methods are one of the most powerful methods that solve hyperbolic problems with the Gibbs oscillations much reduced \cite{Harten,WENO,ShuSIAM}. In \cite{GuoJung}, we developed non-polynomial ENO and WENO methods based on the non-polynomial reconstruction for the finite volume approximations. The original ENO and WENO methods are based on the polynomial interpolation and its convex combination for the final reconstruction at the cell boundaries \cite{ShuSIAM}. Thus once the number of cells that are participating in the interpolation is fixed, the order of accuracy is also fixed. In \cite{GuoJung}, we proposed to use non-polynomial bases such as the radial basis functions (RBFs) to improve the fixed order of accuracy of the final reconstruction. The underlying idea is to use non-polynomial functions as the interpolation basis that contain a free parameter. The value of the introduced free parameter is adaptively determined by the local solutions so that the local accuracy is improved and higher order of accuracy than the polynomial order can be achieved. This improvement was made by making the leading error term vanish by the free parameter or at least vanish to a certain order. If the leading error term vanishes, the accuracy is improved particularly when the solution is smooth. In \cite{GuoJung}, it was shown that the ENO/WENO methods with the non-polynomial basis functions improve the accuracy significantly if the problem is smooth. This paper is a direct extension of our previous work for the finite volume ENO/WENO methods to the finite difference ENO/WENO methods. The formulation of the RBF finite difference ENO/WENO methods is basically same as the finite volume methods, while it is more easier to implement to  higher-dimensional problems. In this paper we show that the RBF finite difference ENO/WENO methods work similarly as the RBF finite volume methods and provide supporting numerical examples. 

For the non-polynomial basis, we use the multi-quadric (MQ) RBFs. We note that the hybrid of the RBF method and the WENO method is not new. For example, in \cite{RBF-ENO} a polyharmonic spline was used to solve PDEs on the unstructured grid with the WENO method. The RBF WENO method in \cite{RBF-ENO} took advantages of both the RBF and WENO methods. That is, the hybrid method utilizes the meshless feature of the RBFs, which is beneficial for the unstructured grid. It also took the advantage of the non-oscillatory feature of the WENO method, which efficiently handles discontinuous solutions. However, as the polyharmonic spline was used as the RBF basis, the shape parameter was not free but fixed. In this paper, we are more interested in improving the accuracy rather than the efficiency dealing with the unstructured mesh. In order to improve local accuracy, we adopt infinitely smooth RBFs such as the MQ-RBFs, which contain the free parameter, so-called the shape parameter. By exploiting the shape parameter, we improve the local accuracy. It would be an interesting project to develop a method which utilizes both the meshless feature and the shape parameter of the RBFs with the WENO method on the unstructured grid. We only focus on the accuracy in this paper. 

We restrict our discussion to RBFs having only one shape parameter although it would be possible to have multiple shape parameters. One of the advantages of using the shape parameter is that the polynomial interpolation is a limit case of the RBF interpolation. That is, the RBF interpolation becomes equivalent to the polynomial interpolation if the shape parameters vanish \cite{Larson,Yoon2}. In this paper, we follow the similar procedure in \cite{GuoJung} for the finite difference ENO/WENO reconstruction and solve several hyperbolic problems. 
Since the shape parameters are optimized using all the given information within the given stencil, the method would be oscillatory if the stencil contains discontinuities. To prevent such oscillations, we adopt the monotone polynomial interpolation by measuring the local maxima. To switch back to the regular ENO/WENO method, we only need to have the shape parameter vanish, which makes it easy to implement the RBF-ENO/WENO method in the existing ENO/WENO code. 

The paper is composed of the following sections. In Section 2, we briefly explain the RBF interpolation and the optimization of the shape parameter.  In Section 3, we explain the regular ENO/WENO finite difference method followed by Section 4 where we compare the polynomial interpolation with the RBF interpolation. The interpolation coefficients for $k = 2$ and $k = 3$ are provided. In Section 5, we provide the numerical examples. For the numerical experiments for discontinuous problems, we first use the monotone polynomial interpolation method for the prevention of oscillations. Then we present the numerical examples for both linear and nonlinear problems and for both scalar and system problems in 1D and 2D. In Section 6, we provide a brief conclusion and our future research. 

\section{RBF Interpolation and optimization of the shape parameter}
We briefly explain the RBF interpolation in one space dimension. Suppose that for a domain $\Omega \subset \mathbb{R}$, a data set $\{ (x_i, u_i) \}_{i=1}^N$ is given where the centers $x_i \in \Omega$ are the coordinates and $u_i$ are the function values of the unknown function $u(x)$ at $x = x_i$.  The RBF interpolation is given by a linear combination of the radial basis function, $\phi: \Omega \rightarrow \mathbb{R}$. The kernel $\phi$ at $x = x_i$ is a function of the {\it radial} distance between $x$ and the center $x_i$ and the shape parameters $\varepsilon_i$, i.e. $\phi= \phi (|| x - x_i ||,\varepsilon_i)$ where $|| x - x_i || $ is the distance between $x$ and $x_i$. The Euclidean norm is usually used for $|| x - x_i ||$.  Then the RBF interpolation of $u(x)$ based on $N$ data points, $I_{N}^R u(x)$, is given by 
\begin{equation}
 I_{N}^R u(x)=\sum_{i=1}^N  \lambda_i\phi(||x - x_i ||, \varepsilon_i),
 \label{equation1}
\end{equation}
where $\lambda_i$ are the expansion coefficients to be determined. Here we note that one can add a low order polynomial term in the right hand side of the above equation to make the interpolation consistent up to a certain degree. In this paper, however, we do not include the low order polynomial term in the RBF interpolation.  Using the interpolation condition $I^R_Nu(x_i) = u_i, i = 1, \cdots, N$, the expansion coefficients are determined by the linear system
\begin{equation}
 A\cdot\vec{\lambda}=\vec{u},
\label{equation2}
\end{equation}
where $\vec{\lambda}=[\lambda_1,\dots,\lambda_N]^T$, $\vec{u}=[u_1,\dots,u_N]^T$ and $A_{ij}=\phi(||x_i- x_j ||, \varepsilon_j)$ is the $(ij)$ element of the interpolation matrix $A$. 
If $\phi(r)$ is one of the piecewise smooth RBFs with no shape parameter involved, we will obtain the RBF interpolation explicitly after we solve the linear system. However, if $\phi(r)$ is one of the infinitely smooth RBFs, we still need to determine or pre-assign the value of the shape parameter. One can have shape parameters fixed globally as a constant. Or one can optimize them so that the accuracy of the interpolation is improved. There are various kinds of RBFs \cite{Buhman}. Among those, we use the multi-quadric (MQ) RBF given by 
\begin{equation}
   \phi (||x - x_i||, \varepsilon_i) = \sqrt{1+ \varepsilon_i^2(x - x_i)^2}. 
   \label{MQRBF}
\end{equation}
One could use a different form of the MQ-RBF such as $  \phi (||x - x_i||, \varepsilon_i) = \sqrt{(x - x_i)^2 + \varepsilon_i^2}$, but this form yields a more complicated algebraic expression for the RBF-ENO/WENO formulation than \eqref{MQRBF}.  


In \cite{Buhman}, it has been discussed that the RBF interpolation may yield spectral convergence. In practice, such a fast convergence depends on several conditions such as the optimization of the shape parameters and the data structure. In this paper, we focus on the optimization of the shape parameter to enhance the overall accuracy. This optimization procedure is the key element of the RBF-ENO/WENO formulation. To explain the optimization procedure below, we first assume that the center set $\{x_i\}$ forms a uniform grid and consider the case of $N = 2$. In the following section, similar analysis will be performed for the  RBF-ENO/WENO formulation. It is straightforward to extend the same idea to larger values of $N$. 

Let $x_1 = 0$, $x_2 = \Delta x$ and $\Delta x = x_2 - x_1$ with $N = 2$. Let $u_i$ denote the function value of $u(x)$ at $x = x_i$. For example, $u_{1\over 2} = u(x_{1\over 2}) = u({{\Delta x} \over 2})$. Also we assume that all the shape parameters are variable but same, i.e. $\varepsilon_1 = \varepsilon_2 = \epsilon$. The value of $\epsilon$ can be either real or complex. Define the error function $E(x) = u(x) - I^R_N (x), x \in [x_1, x_2]$. The interpolation based on the Lagrange polynomials yields the second order convergence in $\Delta x$, i.e. $|E(x)| \sim O(\Delta x^2)$. The RBF interpolation, $I^R_N(x)$ is then given by from \eqref{equation1} 
$$
    I^R_N (x) = \lambda_1 \sqrt{1+ \varepsilon^2 x^2} + \lambda_2\sqrt{1+ \varepsilon^2(x -\Delta x)^2},  
$$
where $\lambda_1$ and $\lambda_2$ are to be computed according to \eqref{equation2}. 

We evaluate the error function $E(x)$ at $x = x_{1\over 2}$. With the fixed $\varepsilon$ and the obtained $\lambda_1$ and $\lambda_2$ from \eqref{equation2}, the Taylor series of $E(x_{1\over 2})$  in terms of $\Delta x$ is given by
\begin{equation}
   E(x_{1\over 2}) = {1\over 8}\left( \varepsilon^2 u_{1\over 2} - u''_{1\over 2} \right) \Delta x^2 + O(\Delta x^4), 
\label{error1}
\end{equation}
where the superscript $''$ denotes the second derivative in $x$. The Taylor series of $E(x)$ around $x = x_1$ yields a  similar result but the cubic order term appears as below
\begin{equation}
  E(x_{1\over 2}) =   {1\over 8}\left( \varepsilon^2 u_{1} - u''_{1} \right) \Delta x^2 
  + {1\over {16}} \left( \varepsilon^2 u'_{1} - u'''_{1}\right) \Delta x^3 + O(\Delta x^4). 
\label{error2}
\end{equation}
From \eqref{error1} and \eqref{error2}, we know that it is possible to obtain the $4$th or $3$rd order accuracy if we choose the shape parameter $\varepsilon$ as 
$$
     \epsilon^2 = {{u''_{1\over 2}}\over{u_{1\over 2}}} \mbox{ for } O(\Delta x^4) \quad  \mbox{ or } \quad \epsilon^2 =  {{u''_{1}}\over{u_{1}}} \mbox{ for } O(\Delta x^3). 
$$
Thus the RBF interpolation can achieve higher convergence than the second order expected by the polynomial interpolation. Or at least we know that there exists $\epsilon^2$ that makes the leading error term vanish and helps the interpolation to achieve higher order accuracy. Note, however, that since $u_1$ and $u_2$ are the only given information about the unknown function $u(x)$, the exact value of $u''(x)$ at $x = x_0$ and $ x_{1\over 2}$ and $u(x)$ at $x = x_{1\over 2}$ are not available.  The key idea for the construction of the RBF-ENO/WENO method is that we use the approximation of those values to a certain order based on the function values at the given cells so that we can still obtain the improved order of accuracy.  Table \ref{table_ep0} shows the adaptation condition of $\epsilon$ that makes the leading error term at $x = x_{1\over 2}$ vanish  when expanded in the Taylor series around $x = x_{1\over 2}$ for different values of $N$ with centers $x_i = (i - 1)\Delta x, i = 1, \cdots, N$. The table shows the leading interpolation errors with different $N$ when the interpolation is obtained as the polynomials ($\epsilon = 0$, {\it Polynomial order}) and RBFs when the adaptation is applied ({\it RBF order}). As shown in the table the error of the RBF interpolation is of $O(\Delta x^{N+2})$ for $N = 2$ because of the symmetry and of $O(\Delta x^{N+1})$ for the rest while the error of the polynomial interpolation is of $O(\Delta x^N)$ for every $N$. 

\begin{table}[h]
\renewcommand{\arraystretch}{2.1}
\caption{Adaptation conditions of $\epsilon$ for various values of $N$.}
\begin{center} \footnotesize
\begin{tabular}{|c|c|c|c|} \hline 
N & Adaptation condition for $\epsilon$ & Polynomial order & RBF order\\ \hline 
2 & $\epsilon^2 h_{\frac{1}{2}}- h''_{\frac{1}{2}} =0$ & $O(\Delta x^2)$ & $O(\Delta x^4)$ \\ \hline 
3 & $3 \epsilon^2 h'_{\frac{1}{2}}+ h'''_{\frac{1}{2}} =0$ & $O(\Delta x^3)$ & $O(\Delta x^4)$\\ \hline 
4 & $9\epsilon^4 h_{\frac{1}{2}}-12 \epsilon^2 h''_{\frac{1}{2}} - h^{(4)}_{\frac{1}{2}}=0$ & $O(\Delta x^4)$ & $O(\Delta x^5)$\\ \hline
5 & $75\epsilon^4 h'_{\frac{1}{2}}+40 \epsilon^2 h'''_{\frac{1}{2}} + h^{(5)}_{\frac{1}{2}}=0$ & $O(\Delta x^5)$ & $O(\Delta x^6)$\\ \hline
\end{tabular}
\end{center} 
\label{table_ep0}
\end{table}

\section{1D Finite difference ENO/WENO method}
Given a uniform grid with $N$ number of points
$$ a = x_{\frac{1}{2}} < x_{\frac{3}{2}} < \cdots < x_{N-\frac{1}{2}} < x_{N+ \frac{1}{2}} = b,$$
for the $i$-th cell $I_i = [x_{i-\frac{1}{2}}, x_{i+\frac{1}{2}}]$, define the cell center $x_i$ and uniform grid spacing $\Delta x$  as
$$x_i = \frac{1}{2} (x_{i-\frac{1}{2}}+ x_{i+\frac{1}{2}}),~  \Delta x = x_{i+\frac{1}{2}} - x_{i-\frac{1}{2}},~  i = 1,2,\cdots,N. $$ 
The semi-discretized form of the 1D hyperbolic equation 
$$ u_t(x,t)+f_x(u(x,t)) = 0 $$
is a system of ordinary differential equations
\begin{eqnarray}
  \frac{du_{i}(t)}{dt} = -\frac{\partial f \left( u(x,t) \right)}{\partial x} \bigg |_{x=x_i}
\label{eq2}
\end{eqnarray} 
where $u_i(t)$ is a numerical approximation to $u(x_i,t)$ and $f$ is the flux function. While the finite volume method reconstructs the cell boundary values of the solution, the finite difference method seeks the reconstruction of the flux function at the cell boundaries. Define the numerical flux function $h(x,t)$ for each cell
$$
f(u(x,t)) = \frac{1}{\Delta x} \int^{x+\frac{\Delta x}{2}}_{x-\frac{\Delta x}{2}} h(\xi,t) d\xi. 
$$
Differentiate both sides and evaluate them at $x=x_i$ to obtain
$$
\frac{\partial f(u(x,t)}{\partial x} \bigg |_{x=x_i} = \frac{1}{\Delta x} \left[ h(x_{i+{1\over 2}},t) - h(x_{i-{1\over 2}},t) \right].
$$
Therefore, \eqref{eq2} becomes 
\begin{eqnarray}
  \frac{du_{i}(t)}{dt} = -\frac{1}{\Delta x} \left[ h_{i+{1\over 2}} - h_{i-{1\over 2}} \right],
\label{eq3}
\end{eqnarray} 
where $h_{i\pm{1\over 2}} = h(x_{i\pm{1\over 2}},t)$.
Given grid values of $f$, $f_i = f(u(x_i,t)) = \frac{1}{\Delta x} \int^{x_{i+{1 \over 2}}}_{x_{i-{1 \over 2}}} h(\xi,t) d\xi$, the high order reconstruction of $\hat{f}_{i\pm{1\over 2}} = h_{i\pm{1\over 2}} + O(\Delta x^m)$ needs to be computed, where $m$ depends on the number of cells we use. 

For the $k$-th order ENO reconstruction, we choose the stencil based on $r$ cells to the left and $s$ cells to the right including $I_i$ such that 
$$
   r + s + 1 = k. 
$$
Define $S_r(i)$ as the stencil composed of those $k$ cells including the cell $I_i$
\begin{equation}
      S_r(i) = \left\{ I_{i-r}, \cdots, I_{i+s} \right\}, ~\quad r=0, \cdots, k-1.
      \label{stencil}
\end{equation}
Define a primitive function $H(x)$ such that 
\begin{equation}
      H(x) = \int^x_{x_{i-r-{1\over 2}}} h(\xi) d\xi 
      \label{primitive}
\end{equation}
where the lower limit in the integral can be any cell boundary \cite{Shu}. By the definition of $H(x)$ in \eqref{primitive}, it is obvious that $H'(x) = h(x)$. Then for $i - r - 1 \leqslant l \leqslant i + s$, $H(x_{l+{1\over 2}})$ is given by the linear sum of cell averages 
$$
   H(x_{l + {1\over 2}}) = \sum_{j = i-r}^l\int^{x_{j+{1\over 2}}}_{x_{j - {1\over 2}}} h(\xi)d\xi = \sum_{j = i-r}^l \Delta x f_j. 
$$
The regular ENO method constructs the polynomial interpolation of $H(x)$ based on $H(x_{l + {1\over 2}}),~ i - r - 1 \leqslant l \leqslant i + s$, while the RBF-ENO method constructs the RBF interpolation of $H(x)$.
Suppose $P(x)$ is some interpolation for $H(x)$ such that 
\begin{equation}
      P(x) =  H(x) + O(\Delta x^{k+1}).
      \label{Px}
\end{equation}
Then $p(x) \equiv P'(x)$ is the function we seek to approximate $h(x)$ where 
\begin{equation}
      p(x) =  h(x) + O(\Delta x^{k}).
      \label{px}
\end{equation}
Then the numerical flux function is given by the linear combination of the given flux value at each point 
$${\hat{f}}^{(r)}_{i+ {1\over 2}} \equiv p(x_{i+{1\over 2}}) = \sum^{k-1}_{j=0} c_{rj} f_{i-r+j}$$ 
where $c_{rj}$ are fixed constants for the ENO method but are variables for the RBF-ENO method.

The WENO method seeks the reconstruction as a convex combination of all possible ENO reconstructions. Suppose the stencil \eqref{stencil} produces $k$ different reconstructions to the value ${\hat{f}}^{(r)}_{i+{1\over 2}}$. Then the WENO reconstruction would take the convex combination of all ${\hat{f}}^{(r)}_{i+ {1\over 2}}$:
\begin{eqnarray}
   {\hat{f}}_{i+{1\over 2}} = \sum^{k-1}_{r=0} w_r {\hat{f}}^{(r)}_{i+ {1\over 2}},
\label{weno_re}
\end{eqnarray}
where 
$$
   w_r = \frac{\alpha_r}{\sum^{k-1}_{s=0} \alpha_s}, r=0,\cdots ,k-1
$$
with 
$$
   \alpha_r = \frac{d_r}{(\epsilon_M+\beta_r)^2}.
$$
Here $d_r$ are the polynomial expansion coefficients and $\epsilon_M>0$ is introduced to avoid the case that the denominator becomes zero usually taken as $\epsilon_M = 10^{-6}$. $\beta_r$ are the ``smoothness indicators" of the stencil $S_r(i)$. In this work, we use the smoothness indicators developed in \cite{WENO}. The usage of different  smoothness indicators found in other WENO variations such as the WENO-Z method \cite{WENOZ}, WENO-M method \cite{WENOM} or WENO-P method \cite{Ha} will be considered in our future work.

\section{1D reconstruction based on RBFs} 
The 1D and 2D finite difference RBF-ENO method is a direct extension of the finite volume RBF-ENO method \cite{GuoJung}. We only replace the cell averages $\bar u_i$ with the nodal flux function $f_i$. First we consider the case of $k=2$, that is, two cells are used for the reconstruction. For this case, three flux values, $f_{i-1}$, $f_{i}$ and $f_{i+1}$ are available. To reconstruct the boundary values of $f_{i+ \frac{1}{2}}$ and $f_{i- \frac{1}{2}}$, we use either $\{ f_{i-1}, f_{i} \}$ or $\{ f_{i}, f_{i+1}\}$. Which flux values should be used is decided by the Newton's divided difference method \cite{Shu}. Assume that we decided to use $\{ f_{i}, f_{i+1}\}$ from the Newton's divided difference. We only show the reconstruction at cell boundary $x = x_{i+ {1\over 2}}$. The reconstruction at $x = x_{i- {1\over 2}}$ can be achieved in the same manner. Based on the definition of $H(x)$ in \eqref{primitive} we have
\begin{eqnarray}
H(x_{i-{1\over 2}}) = - \Delta x \cdot f_{i},  \quad
H(x_{i+{1\over 2}}) = 0, \quad
H(x_{i+{3\over 2}}) = \Delta x \cdot f_{i+1}. \nonumber
\end{eqnarray}
\subsection{Polynomial interpolation for regular ENO}
The polynomial interpolation $P(x)$ is given as 
$$ P(x) = \lambda_1 + \lambda_2 x + \lambda_3 x^2. $$
Let $ {\vec{H}} =  [H_{i- {1\over 2}}, H_{i+ {1\over 2}}, H_{i+ {3\over 2}} ]^T$, $ {\vec{\lambda}} = [\lambda_1, \lambda_2,  \lambda_3]^T$ and the interpolation matrix $A$ be
$$
A = \begin{bmatrix} 
1 & x_{i-{1\over 2}} & x_{i-{1\over 2}}^2
\\ 
1 & x_{i+{1\over 2}} & x_{i+{1\over 2}}^2
\\ 
1 & x_{i+{3\over 2}} & x_{i+{3\over 2}}^2
\end{bmatrix}. 
$$
Then the expansion coefficients $\lambda_i$ are given by solving the linear system $ \vec{H} = A \cdot \vec{\lambda} $.
After taking the first derivative of $P(x)$ and plugging it in $x = x_{i+ {1\over 2}}$, we obtain the numerical flux function ${\hat f}_{i+{1\over 2}}$
\begin{equation}
{\hat f}_{i+ {1\over 2}} \equiv p(x_{i+ {1\over 2}}) = {1\over 2} \cdot f_i+{1\over 2} \cdot f_{i+1}. 
\label{poly_recon}
\end{equation}
Expanding ${\hat f}_{i+ {1\over 2}}$ around $x=x_{i+{1\over 2}}$ in the Taylor series yields
\begin{equation}
{\hat f}_{i+ {1\over 2}} \equiv p(x_{i+{1\over 2}}) = h_{i+{1\over 2}} + {1\over 6}  h''_{i+{1\over 2}} \Delta x^2 + {1\over 120} h^{(4)}_{i+{1\over 2}} {\Delta x^4} + O(\Delta x^6). 
\label{poly_recon_error}
\end{equation}
The first term in the right hand side of \eqref{poly_recon_error} is the exact value of $h(x)$ at $x = x_{i+ {1\over 2}}$. Thus we confirm that \eqref{poly_recon} is a $2$nd order reconstruction of $f$.

\subsection{Multi-Quadratic RBF interpolation}
Now we consider the RBF reconstruction of $f$. The MQ-RBF interpolation $P(x)$ for $H(x)$ is given as
$$ P(x) = \lambda_1 \sqrt{1+\epsilon^2(x-x_{i- {1\over 2}})^2  } + \lambda_2 \sqrt{1+\epsilon^2(x-x_{i+ {1\over 2}})^2  } + \lambda_3 \sqrt{1+\epsilon^2(x-x_{i+ {3\over 2}})^2  }.  $$
Accordingly the interpolation matrix $A$ is given by 
$$
A = \begin{bmatrix} 
1 & \sqrt{\Delta x^2 \epsilon^2+1} & \sqrt{4 \Delta x^2 \epsilon^2+1}
\\ 
\sqrt{\Delta x^2 \epsilon^2+1} & 1 & \sqrt{\Delta x^2 \epsilon^2+1}
\\ 
\sqrt{4 \Delta x^2 \epsilon^2+1} & \sqrt{\Delta  x^2 \epsilon^2+1} & 1
\end{bmatrix}. 
$$
We take the first derivative of $P(x)$ to obtain the numerical flux at $x_{i+{1\over 2}}$
\begin{equation}
{\hat f}_{i+ {1\over 2}} = \left(\frac{\sqrt{4 \epsilon^2 \Delta x^2  +1}+1}{4 \sqrt{\epsilon^2 \Delta x^2 +1} }  \right) \cdot f_i+ \left(\frac{\sqrt{4 \epsilon^2 \Delta x^2 +1}+1}{4 \sqrt{\epsilon^2 \Delta x^2 +1} } \right) \cdot f_{i+1}. 
\label{i+1/2_exact}
\end{equation}
Expanding ${\hat f}_{i+ {1\over 2}}$ around $x=x_{i+{1\over 2}}$ in the Taylor series yields
\begin{eqnarray}
{\hat f}_{i+ {1\over 2}} &=& h_{i+{1\over 2}} + \left({1\over 2}\epsilon^2 h_{i+{1\over 2}}+{1\over 6} h''_{i+{1\over 2}} \right) \Delta x^2 
\nonumber \\
& & + \left( -{9\over 8} \epsilon^4  h_{i+{1\over 2}}  + {1\over 12} \epsilon^2  h''_{i+{1\over 2}} + {1\over 120} h^{(4)}_{i+{1\over 2}}  \right) {\Delta x^4} + O(\Delta x^6). 
\label{p(i+1/2)for_error}
\end{eqnarray}
Thus we observe that \eqref{i+1/2_exact} is at least $2$nd order accurate to $h_{i+{1\over 2}}$.  If we take the value of $\epsilon$ as below
\begin{equation}
\epsilon^2 = -{1\over 3}{{h''_{i+{1\over 2}}}\over{h_{i+{1\over 2}}}},
\label{i+1/2_eps}
\end{equation}
then the $2$nd order error term vanishes and we obtain a $4$th order accurate approximation. 
We notice that the coefficients of $f_i$ and $f_{i+1}$ in \eqref{i+1/2_exact} have complicated forms involving a calculation of square roots. This is common among infinitely smooth RBF interpolations because of the existence of undetermined shape parameters. If we fix the shape parameter $\epsilon$ as a constant, the reconstruction would have simpler forms as the piecewise smooth RBF interpolation, but lose the ability to improve accuracy. In order to simplify the form, we expand the right hand side of \eqref{i+1/2_exact} as below
\begin{eqnarray}
{\hat f}_{i+ {1\over 2}}&=& \left(\frac{1}{2}+\frac{1}{4} \epsilon^2 \Delta x^2-\frac{9}{16} \epsilon^4 \Delta x^4\right) \cdot f_i + \left(\frac{1}{2}+\frac{1}{4} \epsilon^2 \Delta x^2-\frac{9}{16} \epsilon^4 \Delta x^4\right) \cdot f_{i+1} \nonumber \\
  &&+ O(\Delta x^6). 
\label{i+1/2_around_dx}
\end{eqnarray}
Ignoring all the high order terms in \eqref{i+1/2_around_dx} yields
\begin{equation}
{\hat f}_{i+ {1\over 2}} = \left(\frac{1}{2}+\frac{1}{4} \epsilon^2 \Delta x^2\right) \cdot f_i + \left(\frac{1}{2}+\frac{1}{4} \epsilon^2 \Delta x^2\right) \cdot f_{i+1}. 
\label{i+1/2_approx}
\end{equation}
Expanding $f_i$ and $f_{i+1}$ in terms of $h(x)$ we get 
\begin{eqnarray}
{\hat f}_{i+ {1\over 2}}  &=& h_{i+{1\over 2}} + \left({1\over 2}\epsilon^2  h_{i+{1\over 2}}  + {1\over 6}  h''_{i+{1\over 2}} \right) \Delta x^2  
\nonumber \\
&& + \left( {1\over 12} \epsilon^2  h''_{i+{1\over 2}} + {1\over 120} h^{(4)}_{i+{1\over 2}} \right) {\Delta x^4} + O(\Delta x^6). 
\label{p(i+1/2)for_error_approx}
\end{eqnarray}
Thus \eqref{i+1/2_approx} yields the same $4$th order accuracy equivalently to \eqref{i+1/2_exact} when \eqref{i+1/2_eps} is used. 

Now we consider the evaluation of $\epsilon$ in \eqref{i+1/2_eps} to achieve higher order accuracy than the $2$nd order. To explain this, we use the reconstruction at  $x = x_{i+{1\over 2}}$ as an example. The case of $x=x_{i-{1\over 2}}$ can be done in the same manner.
The idea is that although the ENO method chooses either $\{f_i, f_{i+1}\}$ or $\{f_{i}, f_{i-1}\}$ for the final reconstruction, all three flux function values are available when the Newton's divided difference method is applied. Thus we calculate $h_{i+{1\over 2}}$ and $h''_{i+{1\over 2}}$ based on $f_{i-1}, f_i$ and  $f_{i+1}$. 
Construct the Lagrange interpolation of $H(x)$ based on $H_{i-{3\over 2}}, H_{i-{1\over 2}}, H_{i+{1\over 2}},H_{i+{3\over 2}}$ and take the first derivative of $H(x)$ at $x = x_{i+{1\over 2}}$. Then we have 
$$
    h_{i+{1\over 2}} \approx -{1\over 6} f_{i-1} + {5\over 6} f_{i} + {1\over 3} f_{i+1} + O(\Delta x^3). 
$$
The second derivative of $ h_{i+{1\over 2}} $ is approximated by the third derivative of $H(x)$, which is given by
$$
  h''_{i+{1\over 2}} \approx {{f_{i-1} - 2 f_{i} + f_{i+1}}\over{\Delta x^2}}+ O(\Delta x). 
$$
Then by plugging the above approximations of $h_{i+{1\over 2}}$ and $h''_{i+{1\over 2}}$ into \eqref{i+1/2_eps}, we approximate the optimal value of $\epsilon^2$ as below
\begin{equation}
\epsilon^2 \approx  \frac{2}{\Delta x^2}\cdot \frac{-f_{i-1}+2f_{i}-f_{i+1}}{-f_{i-1}+5f_{i}+2 f_{i+1}}.
\label{i+1/2_eps_approx}
\end{equation}
Note that $\epsilon$ can be a complex number because $\epsilon^2$ can be a negative real number. Since $\epsilon^2$ is used instead of $\epsilon$, the final reconstruction is still done with the real operation. 
If we replace $\epsilon^2$ in \eqref{p(i+1/2)for_error} with that in \eqref{i+1/2_eps_approx}, we get the following 
$$ {\hat f}_{i+ {1\over 2}} \equiv p(x_{i+{1\over 2}}) = h_{i+{1\over 2}} + {1\over 12} h^{(3)}_{i+{1\over 2}} \Delta x^3 + O(\Delta x^4). $$
Thus if $\epsilon^2$ is chosen as in \eqref{i+1/2_eps_approx}, the reconstruction \eqref{i+1/2_approx} still yields a  $3$rd order accuracy although we do not completely make the $2$nd order term in \eqref{p(i+1/2)for_error_approx} vanish.

\subsection{Reconstruction coefficients}
Tables \ref{table_k2} and \ref{table_k3} show the reconstruction coefficients for $k = 2$ and $k = 3$, respectively for the MQ RBF-ENO reconstruction. These are the same as for the finite volume method provided in \cite{GuoJung}. The difference is that the flux function is used for the calculation of the optimal value of $\epsilon^2$, while the cell average of the solution is used for the finite volume method. As we see in these tables, the RBF-ENO reconstruction is a perturbed reconstruction of the polynomial-based ENO interpolation. We also observe that the coefficients for $k = 2$ are not consistent while those for $k = 3$ are consistent. Here note that when the shape parameter vanishes, i.e. when $\epsilon \rightarrow 0$, the MQ-RBF reconstruction is reduced into the polynomial reconstruction. For this equivalence known as the polynomial limit \cite{Larson,Yoon2}, one can easily switch the RBF-ENO to the regular ENO and it is easy to convert the existing ENO code into the RBF-ENO code. 

Table \ref{table_ep} shows the leading error terms we want to make vanish as in \eqref{p(i+1/2)for_error}  to achieve (k+1)th order of convergence  when $k$ cells are used. The adaptation condition becomes complicated for the approximation of  $\epsilon$ when the value of $k$ becomes large. We leave efficient evaluation methods of $\epsilon$ for larger values of $k$ for our future work. 

\begin{table}[h]
\renewcommand{\arraystretch}{2.1}
\caption{Left: The polynomial reconstruction coefficients. Right: The MQ-RBF reconstruction coefficients  $c_{rj}$ in \eqref{i+1/2_approx}. $\eta =  \epsilon^2 \Delta x^2$.}
\begin{center} \footnotesize
\begin{tabular}{|c|c|c|c|} \hline  
k & r & j=0 & j=1\\ \hline 
  & -1 &  $\frac{3}{2}$ &  $-\frac{1}{2}$\\  \cline{2-4}
2 & 0 &  $\frac{1}{2} $ & $\frac{1}{2} $ \\  \cline{2-4} 
  & 1 & $-\frac{1}{2}$ &  $\frac{3}{2}$ \\   \cline{1-4} 
\multicolumn{4}{|c|}{$\eta = 0$} \\
\hline
\end{tabular}
\hskip .2in
\begin{tabular}{|c|c|c|c|} \hline  
k & r & j=0 & j=1\\ \hline 
  & -1 & ${3\over 2} - {3\over 2} \eta$ & $- {1\over 2} + {1\over 2} \eta$\\  \cline{2-4}
2 & 0 &  ${1\over 2} + {1\over 4} \eta$ & $ {1\over 2} + {1\over 4} \eta$ \\  \cline{2-4} 
  & 1 & $- {1\over 2} + {1\over 2} \eta$ &  ${3\over 2} - {3\over 2} \eta$ \\
\cline{1-4} 
\multicolumn{4}{|c|}{$\eta = \epsilon^2 \Delta x^2  = \frac{ 2(-f_{i-1}+2f_{i}-f_{i+1})}{-f_{i-1}+5f_{i}+2f_{i+1}}$} \\
\hline
\end{tabular}
\end{center} 
\label{table_k2}


\renewcommand{\arraystretch}{2.1}
\caption{Left: The polynomial reconstruction. Right: The MQ-RBF reconstruction coefficients  $c_{rj}$, $\eta =  \epsilon^2 \Delta x^2$.}
\begin{center} \footnotesize
\begin{tabular}{|c|c|c|c|c|} \hline  
k & r & j=0 & j=1 & j=2\\ \hline 
  & -1 & $\frac{6}{11}$ & $-\frac{7}{6}$ & $\frac{1}{3}$ \\  \cline{2-5}
3 & 0 &  $\frac{1}{3}$ & $\frac{5}{6}$ & $-\frac{1}{6}$ \\  \cline{2-5} 
  & 1 & $-\frac{1}{6}$ & $\frac{5}{6}$ & $\frac{1}{3}$\\  \cline{2-5}
  & 2 & $\frac{1}{3}$ & $-\frac{7}{6}$ & $\frac{11}{6}$ \\  \cline{1-5} 
\multicolumn{5}{|c|}{$\eta = 0$} \\
\hline
\end{tabular} 
\hskip .2in
\begin{tabular}{|c|c|c|c|c|} \hline  
k & r & j=0 & j=1 & j=2\\ \hline
  & -1 & $\frac{6}{11} - \frac{9}{2} \eta$ & $-\frac{7}{6} + 6 \eta$ & $\frac{1}{3} - \frac{3}{2} \eta$ \\  \cline{2-5}
3 & 0 &  $\frac{1}{3} + \frac{5}{6} \eta$ & $\frac{5}{6} - \frac{2}{3} \eta$ & $-\frac{1}{6} - \frac{1}{6} \eta$ \\  \cline{2-5} 
  & 1 & $-\frac{1}{6} - \frac{1}{6} \eta$ & $\frac{5}{6} - \frac{2}{3} \eta$ & $\frac{1}{3} + \frac{5}{6} \eta$\\  \cline{2-5}
  & 2 & $\frac{1}{3} - \frac{3}{2} \eta$ & $-\frac{7}{6} + 6 \eta$ & $\frac{11}{6} - \frac{9}{2} \eta$ \\  \cline{1-5} 
\multicolumn{5}{|c|}{$\eta =\epsilon^2 \Delta x^2=  {{f_{i-1} - 3f_i + 3 f_{i+1} - f_{i+2}}\over{f_{i-1} - 15 {\bar{v}}_i + 15 f_{i+1} - f_{i+2} } } $} \\
\hline
\end{tabular} 

\end{center} 
\label{table_k3}
\end{table}

\begin{table}[h]
\renewcommand{\arraystretch}{2.1}
\caption{Adaptation conditions of $\epsilon$ for various values of $k$. $h_j = h(x_{i + {1 \over 2}})$ or $h(x_{i - {1 \over 2}}).$  }
\begin{center} \footnotesize
\begin{tabular}{|c|c|} \hline 
k & coefficient in the leading error term $O( \Delta x^k)$ \\ \hline 
2 & $\epsilon^2 h_j+{1\over 3} h''_j =0$ \\ \hline 
3 & $\epsilon^2 h'_j+{1\over 12} h'''_j =0$\\ \hline 
4 & $15\epsilon^4 h_j+8 \epsilon^2 h''_j + {1 \over 5} h^{(4)}_j=0$\\ \hline
5 & $1275\epsilon^4 h'_j+185 \epsilon^2 h'''_j + {7 \over 3} h^{(5)}_j=0$\\ \hline
\end{tabular}
\end{center} 
\label{table_ep}
\end{table}

\section{Numerical experiments}
For the numerical example, we compare the developed method with the original ENO/WENO method denoted by the ENO/WENO-JS developed by Jiang and Shu \cite{WENO}. For the WENO method we use the smoothness indicators developed in \cite{WENO} and denote our proposed RBF-WENO method as the RBF-WENO-JS method.  One can adopt other developments rather than the original WENO method for the RBF-WENO method such as the WENO-Z method \cite{WENOZ}, which will be investigated in our future research. For the numerical experiments, we use the Lax-Friedrich flux scheme with the $3$rd order TVD Runge-Kutta method \cite{RK3} for the time integration. 

\subsection{Polynomial limit and non-oscillatory reconstruction}
The ENO method uses adaptive stencil to avoid the reconstruction across the possible discontinuity. The RBF-ENO method, however, as the WENO reconstruction, utilizes all available flux functions to optimize the shape parameter within the stencil $S_r(i)$. Thus the final reconstruction can be oscillatory. To prevent the oscillation and recover the ENO property, one can easily switch the RBF-ENO reconstruction to the ENO reconstruction by using the polynomial limit of RBFs. That is, if the shape parameter $\epsilon$ vanishes, the RBF reconstruction becomes equivalent to the polynomial reconstruction and becomes the regular ENO reconstruction. In \cite{GuoJung} the monotone polynomial method was proposed to use the polynomial limit whenever it is necessary. For the finite difference RBF-ENO method, we follow the same procedure. The difference is that for the local extrema we use the flux function $f_i$ given at each $x_i$ instead of the solution $u_i$. 

The optimized value of $\epsilon^2$ in \eqref{i+1/2_eps} involves the second derivative of $h$, e.g. $h''_{i+{1\over 2}}$. For $k = 2$, there are three cells involved and only one $h''$ is to be computed while for $k = 3$, three possible values of $h''$ are to be computed. For the illustration consider the case of $k = 2$. For $k = 3$, the same procedure is straightforwardly applied for each ENO block.  For $k = 2$, 
%
%
the whole interval of $x$ in the stencil is $ x_{i-{3\over 2}} \le x \le x_{i+{3\over 2}}$.  Without loss of generality, let $x_{i-{3\over 2}} = 0$ and $x_{i+{3\over 2}} = 3\Delta x$.  In this interval,  $h(x)$ is approximated by the second order polynomial and its local extrema exists at $x = x_p$
%
%
\begin{eqnarray}
x_p = {{-2f_{i-1} +3f_{i}-f_{i+1} }\over{-f_{i-1} +2f_{i}-f_{i+1}}} \Delta x.
\label{root_value}
\end{eqnarray}
If $x_p$ exists inside the given stencil, $p(x)$ is not a monotone function.
We use the polynomial limit if $x_p$ exists inside the interval of $S_r(i)$. Thus the polynomial limit condition is given by 
\begin{eqnarray}
  \epsilon^2 = \left\{ \begin{array}{ll} 0  & \mbox{ if }  {0} < {{-2f_{i-1} +3f_{i}-f_{i+1} }\over{-f_{i-1} +2f_{i}-f_{i+1}}}  < {3} \\
   -{1\over 3}{{  h''_{i+{1\over 2}}}\over{h_{i+{1\over 2}}}} & \mbox{ otherwise } \end{array} \right. .
\label{adaptation}
\end{eqnarray}
\subsection{1D Numerical examples} 

\subsubsection{Advection equation I: smooth initial condition}
First we consider the linear scalar equation in $x \in [-1,1]$
%
%
\begin{eqnarray}
u_t + u_x = 0, \quad  t > 0,
\label{example1}
\end{eqnarray}
with the initial condition 
\begin{eqnarray}
u (x,0) = u_0(x) = \sin(\pi x)
\label{example1_ic}
\end{eqnarray}
and the periodic boundary condition. The CFL condition is given by $\Delta t \le C \Delta x$ with $C = 0.1$. Tables \ref{table0} and \ref{table1} show the $L_1, L_2$ and $L_\infty$ errors for each method at the final time $T = 0.5$ with $k = 2$ and $k = 3$, respectively. 
As shown in the tables, for $k = 2$, the RBF-ENO method has almost $3$rd order convergence while the regular ENO method yields $2$nd order of accuracy or less. The RBF-ENO is also better performed than the WENO-JS in terms of accuracy while the rates of convergence are similar. The RBF-WENO-JS with $k = 2$ is also better than WENO-JS both in accuracy and convergence. We have similar results for $k = 3$. For $k = 3$ the RBF-ENO method yields convergence higher than $3$rd order and the RBF-WENO-JS method higher than $5$th order while the regular ENO method is about $3$rd order and the WENO-JS is about $5$th order accurate. 

\begin{table}[h]
\caption{$L_1, L_2$ and $L_\infty$ errors for the advection equation  \eqref{example1} with continuous initial condition \eqref{example1_ic}. $T = 0.5$. $k = 2$. }
\begin{center} \footnotesize
\renewcommand{\arraystretch}{1}
\begin{tabular}{|c|c|c|c|c|c|c|c|} 
\hline  
Method & N & $L_1$ error & $L_1$ order & $L_2$ error & $L_2$ order & $L_\infty$ error & $L_\infty$ order\\ 
\hline 
       & 10 & 1.11E-1 &   --   & 1.40E-1 &    --  & 2.22E-1 & --     \\  
       & 20 & 4.61E-2 & 1.2689 & 5.32E-2 & 1.3994 & 9.43E-2 & 1.2375 \\  
ENO    & 40 & 1.37E-2 & 1.7458 & 1.79E-2 & 1.5749 & 4.03E-2 & 1.2236 \\
k = 2  & 80 & 3.81E-3 & 1.8530 & 5.69E-3 & 1.6492 & 1.68E-2 & 1.2614 \\ 
       & 160& 1.02E-3 & 1.8901 & 1.80E-3 & 1.6599 & 6.91E-3 & 1.2849 \\ 
       & 320& 2.71E-4 & 1.9224 & 5.69E-4 & 1.6613 & 2.81E-3 & 1.2964 \\  
\hline
           & 10 & 1.79E-2 &    --  & 2.35E-2 &    --  & 4.24E-2 &    --  \\  
           & 20 & 2.48E-3 & 2.8540 & 2.65E-3 & 3.1473 & 3.62E-3 & 3.5495 \\  
RBF-ENO    & 40 & 3.17E-4 & 2.9673 & 3.43E-4 & 2.9488 & 4.79E-4 & 2.9207 \\
k = 2      & 80 & 4.05E-5 & 2.9704 & 4.42E-5 & 2.9592 & 6.22E-5 & 2.9434 \\ 
           & 160& 5.17E-6 & 2.9702 & 5.60E-6 & 2.9810 & 7.97E-6 & 2.9652 \\ 
           & 320& 6.51E-7 & 2.9897 & 7.05E-7 & 2.9892 & 1.00E-6 & 2.9894 \\
\hline
           & 10 & 9.13E-2 &   --   & 1.10E-1 &    --  & 1.74E-1 &   --   \\  
           & 20 & 2.92E-2 & 1.6453 & 3.26E-2 & 1.7588 & 5.51E-2 & 1.6547 \\  
WENO-JS    & 40 & 4.81E-3 & 2.5998 & 6.33E-3 & 2.3645 & 1.38E-2 & 2.0037 \\
k = 2      & 80 & 6.42E-4 & 2.9066 & 9.43E-4 & 2.7470 & 2.61E-3 & 2.3979 \\ 
           & 160& 7.79E-5 & 3.0440 & 1.23E-4 & 2.9427 & 3.96E-4 & 2.7183 \\ 
           & 320& 9.54E-6 & 3.0291 & 1.53E-5 & 3.0035 & 5.25E-5 & 2.9181 \\
\hline
           & 10 & 1.88E-2 &    --  & 2.38E-2 &    --  & 4.06E-2 &    --  \\  
           & 20 & 2.60E-3 & 2.8561 & 2.73E-3 & 3.1235 & 4.13E-3 & 3.2995 \\  
RBF-WENO-JS   & 40 & 3.25E-4 & 2.9969 & 3.57E-4 & 2.9367 & 5.36E-4 & 2.9442 \\
k = 2      & 80 & 4.05E-5 & 3.0064 & 4.50E-5 & 2.9871 & 6.75E-5 & 2.9896 \\ 
           & 160& 5.09E-6 & 2.9910 & 5.63E-6 & 2.9983 & 8.32E-6 & 3.0199 \\ 
           & 320& 6.40E-7 & 2.9927 & 7.04E-7 & 3.0000 & 1.00E-6 & 3.0539 \\
\hline
\end{tabular}
\end{center} 
\label{table0}
\bigskip
\caption{$L_1, L_2$ and $L_\infty$ errors for the advection equation  \eqref{example1} with continuous initial condition \eqref{example1_ic}. $T = 0.5$. $k = 3$.}
\begin{center} \footnotesize
\renewcommand{\arraystretch}{1}
\begin{tabular}{|c|c|c|c|c|c|c|c|} 
\hline  
Method & N & $L_1$ error & $L_1$ order & $L_2$ error & $L_2$ order & $L_\infty$ error & $L_\infty$ order\\ 
\hline 
       & 10 & 2.32E-2 &   --   & 2.54E-2 &    --  & 3.65E-2 & --     \\  
       & 20 & 2.79E-3 & 3.0570 & 3.02E-3 & 3.0727 & 4.47E-3 & 3.0314 \\  
ENO    & 40 & 3.36E-4 & 3.0519 & 3.69E-4 & 3.0360 & 5.47E-4 & 3.0293 \\
k = 3  & 80 & 4.12E-5 & 3.0284 & 4.55E-5 & 3.0189 & 6.76E-5 & 3.0167 \\ 
       & 160& 5.10E-6 & 3.0146 & 5.65E-6 & 3.0099 & 8.53E-6 & 2.9878 \\ 
       & 320& 6.34E-7 & 3.0074 & 7.03E-7 & 3.0051 & 1.06E-6 & 3.0082 \\  
\hline
           & 10 & 1.89E-2 &   --   & 2.22E-2 &    --  & 3.45E-2 &   --   \\  
           & 20 & 2.15E-3 & 3.1332 & 2.59E-3 & 3.1022 & 4.47E-3 & 2.9462 \\  
RBF-ENO    & 40 & 1.51E-4 & 3.8321 & 2.18E-4 & 3.5731 & 5.14E-4 & 3.1218 \\
k = 3      & 80 & 8.81E-6 & 4.1023 & 1.57E-5 & 3.7941 & 5.10E-5 & 3.3319 \\ 
           & 160& 4.83E-7 & 4.1891 & 1.06E-6 & 3.8854 & 4.60E-6 & 3.4706 \\ 
           & 320& 2.78E-8 & 4.1205 & 7.30E-8 & 3.8618 & 4.25E-7 & 3.4372 \\  
\hline
           & 10 & 9.74E-3 &   --   & 1.13E-2 &    --  & 1.63E-2 &   --   \\  
           & 20 & 4.01E-4 & 4.6032 & 4.63E-4 & 4.6147 & 7.82E-4 & 4.3822 \\  
WENO-JS    & 40 & 1.17E-5 & 5.0966 & 1.37E-5 & 5.0804 & 2.45E-5 & 4.9965 \\
k = 3      & 80 & 3.55E-7 & 5.0452 & 4.10E-7 & 5.0617 & 7.55E-7 & 5.0196 \\ 
           & 160& 1.12E-8 & 4.9825 & 1.27E-8 & 5.0100 & 2.33E-8 & 5.0180 \\ 
           & 320& 3.80E-10& 4.8849 & 4.25E-10& 4.9022 & 6.94E-10& 5.0692 \\ 
\hline     
		   & 10 & 6.53E-3 &   --   & 7.90E-3 &    --  & 1.16E-2 &   --   \\  
           & 20 & 1.22E-5 & 5.7476 & 1.50E-4 & 5.7155 & 2.72E-4 & 5.4183 \\  
RBF-WENO-JS   & 40 & 2.58E-6 & 5.5559 & 3.24E-6 & 5.5380 & 8.04E-6 & 5.0796 \\
k = 3      & 80 & 7.51E-8 & 5.1050 & 8.72E-8 & 5.2129 & 2.17E-7 & 5.2082 \\ 
           & 160& 2.35E-9 & 4.9964 & 2.63E-9 & 5.0497 & 5.58E-9 & 5.2833 \\ 
           & 320& 7.39E-11& 4.9920 & 8.22E-11& 5.0014 & 1.46E-10& 5.2581 \\  
\hline
\end{tabular}
\end{center} 
\label{table1}
\end{table}

\subsubsection{Advection equation II: non-smooth initial condition}
We consider the same advection equation \eqref{example1} but with the discontinuous initial condition 
\begin{eqnarray}
   u(x,0) = -\mbox{sgn}(x),
   \label{example2} 
\end{eqnarray}
and the boundary condition $u(-1,t) = 1, t>0$. Here $\mbox{sgn(x)}$ is the sign function of $x$. 
Figure \ref{figure2} shows solutions at $T = 0.5$ by each method with $N = 200$. The top two figures show the solutions with $k = 2$ and the bottom two figures with $k = 3$. As shown in the figures, the Gibbs oscillations are not significant and the RBF-ENO and RBF-WENO-JS methods yield almost non-oscillatory solutions. For $k = 2$, the RBF-ENO and RBF-WENO-JS solutions are similar and have sharper solution profiles than the regular ENO and WENO methods. For $k = 3$. The RBF-WENO-JS solution is slightly sharper than the WENO-JS solution. The RBF-ENO solution is sharper than the ENO solution but smoother than the RBF-WENO-JS or WENO-JS solutions. 

%
%
%
%
%
%

\subsubsection{Burger's equation}
We consider the Burgers' equation for $x \in [-1,1]$ 
\begin{eqnarray}
   u_t + \left(\frac{1}{2} u^2 \right)_x = 0,  \quad t >0, 
   \label{example3}
\end{eqnarray}
with the initial condition 
\begin{eqnarray}
        u(x,0) = -\sin(\pi x). 
        \label{example3_example}
\end{eqnarray}
Tables \ref{table2} and \ref{table3} show various errors for each case for $k = 2$ and $k = 3$, respectively. For $k = 2$, the RBF-ENO and RBF-WENO-JS methods yield much better accuracy than the ENO or WENO-JS methods, e.g. $L_\infty$ errors. For $k = 3$, the RBF-WENO-JS method yields convergence higher than $5$th order while the WENO-JS is about $5$th order accurate. The RBF-ENO solution yields convergence higher than $3$rd order while the regular ENO solution is about $3$rd order or less accurate. 
\begin{table}[h]
\caption{$L_1, L_2$ and $L_\infty$ errors for the Burgers' equation  \eqref{example3} with continuous initial condition \eqref{example3_example}. $T = 0.2$. $k = 2$. }
\begin{center} \footnotesize
\renewcommand{\arraystretch}{1}
\begin{tabular}{|c|c|c|c|c|c|c|c|} 
\hline  
Method & N & $L_1$ error & $L_1$ order & $L_2$ error & $L_2$ order & $L_\infty$ error & $L_\infty$ order \\
\hline
           & 10 & 5.38E-2 &   --   & 8.24E-2 &    --  & 1.71E-1 &    --  \\  
           & 20 & 1.98E-2 & 1.4464 & 2.81E-2 & 1.5521 & 5.74E-2 & 1.5716 \\  
ENO        & 40 & 5.40E-3 & 1.8715 & 8.28E-3 & 1.7628 & 2.12E-2 & 1.4355 \\
k = 2      & 80 & 1.31E-3 & 2.0396 & 2.52E-3 & 1.7183 & 8.85E-3 & 1.2619 \\ 
           & 160& 3.85E-4 & 1.7690 & 8.13E-3 & 1.6309 & 3.64E-3 & 1.2816 \\ 
           & 320& 1.07E-4 & 1.8475 & 2.59E-4 & 1.6519 & 1.49E-3 & 1.2928 \\  
\hline
           & 10 & 4.13E-2 &   --   & 4.52E-2 &    --  & 1.16E-1 &    --  \\  
           & 20 & 7.41E-3 & 2.4795 & 1.53E-2 & 1.5817 & 4.63E-2 & 1.3254 \\  
RBF-ENO    & 40 & 1.08E-3 & 2.7836 & 3.01E-3 & 2.3418 & 1.18E-2 & 1.9714 \\
k = 2      & 80 & 1.38E-4 & 2.9656 & 4.22E-4 & 2.8362 & 1.81E-3 & 2.7073 \\ 
           & 160& 1.55E-5 & 3.1548 & 4.89E-5 & 3.1107 & 2.40E-4 & 2.9135 \\ 
           & 320& 1.75E-6 & 3.1415 & 5.54E-6 & 3.1396 & 2.67E-5 & 3.1669 \\
\hline
           & 10 & 4.61E-2 &   --   & 7.23E-2 &    --  & 1.55E-1 &    --   \\  
           & 20 & 1.27E-2 & 1.8592 & 1.96E-2 & 1.8792 & 4.43E-2 & 1.8056  \\  
WENO-JS    & 40 & 1.97E-3 & 2.6901 & 3.13E-3 & 2.6493 & 8.18E-3 & 2.4392  \\
k = 2      & 80 & 2.99E-4 & 2.7199 & 5.22E-4 & 2.5851 & 1.40E-3 & 2.5473  \\ 
           & 160& 4.32E-5 & 2.7893 & 7.89E-5 & 2.7252 & 2.49E-4 & 2.4893  \\ 
           & 320& 5.75E-6 & 2.9100 & 1.07E-5 & 2.8769 & 3.68E-5 & 2.7572  \\
\hline 
           & 10 & 2.50E-2 &   --   & 3.69E-2 &    --  & 7.02E-2 &    --  \\  
           & 20 & 6.64E-3 & 1.9158 & 1.57E-2 & 1.2310 & 4.88E-2 & 0.5227 \\  
RBF-WENO-JS   & 40 & 1.05E-3 & 2.6662 & 3.06E-3 & 2.3610 & 1.28E-2 & 1.9376 \\
k = 2      & 80 & 1.30E-4 & 3.0062 & 4.03E-4 & 2.9230 & 1.70E-3 & 2.9098 \\ 
           & 160& 1.50E-5 & 3.1181 & 4.62E-5 & 3.1254 & 2.15E-4 & 2.9772 \\ 
           & 320& 1.73E-6 & 3.1128 & 5.29E-6 & 3.1286 & 2.45E-5 & 3.1374 \\
\hline 
\end{tabular}
\end{center} 
\label{table2}
\bigskip
\caption{$L_1, L_2$ and $L_\infty$ errors for the Burgers' equation  \eqref{example3} with continuous initial condition \eqref{example3_example}. $T = 0.2$. $k = 3$.}
\begin{center} \footnotesize
\renewcommand{\arraystretch}{1}
\begin{tabular}{|c|c|c|c|c|c|c|c|} 
\hline
Method & N & $L_1$ error & $L_1$ order & $L_2$ error & $L_2$ order & $L_\infty$ error & $L_\infty$ order \\
\hline
           & 10 & 2.00E-2 &   --   & 3.63E-2 &    --  & 7.95E-2 &  --     \\  
           & 20 & 4.89E-3 & 2.0328 & 7.77E-3 & 2.2243 & 1.57E-2 & 2.3449  \\  
ENO        & 40 & 8.68E-4 & 2.4953 & 1.50E-3 & 2.3751 & 4.93E-3 & 1.6676  \\
k = 3      & 80 & 1.25E-4 & 2.7961 & 2.21E-4 & 2.7622 & 6.38E-4 & 2.9499  \\ 
           & 160& 1.84E-5 & 2.7675 & 3.44E-5 & 2.6816 & 1.32E-4 & 2.2697  \\ 
           & 320& 2.84E-6 & 2.6918 & 5.68E-6 & 2.5991 & 2.84E-5 & 2.2184  \\  
\hline        
           & 10 & 1.87E-2 &   --   & 3.48E-2 &    --  & 7.70E-2 &   --   \\  
           & 20 & 4.05E-3 & 2.2024 & 6.22E-3 & 2.4851 & 1.44E-2 & 2.4187 \\  
RBF-ENO    & 40 & 6.21E-4 & 2.7069 & 1.17E-3 & 2.4052 & 3.55E-3 & 2.0213 \\
k = 3      & 80 & 1.22E-4 & 2.3453 & 3.00E-4 & 1.9658 & 1.53E-3 & 1.2130 \\ 
           & 160& 1.03E-5 & 3.5705 & 3.00E-5 & 3.3226 & 1.91E-4 & 3.0040 \\ 
           & 320& 8.23E-7 & 3.6436 & 2.76E-6 & 3.4444 & 2.09E-5 & 3.1898 \\ 
\hline
           & 10 & 1.28E-2 &   --   & 2.21E-2 &    --  & 4.75E-2 &   --    \\  
           & 20 & 2.67E-3 & 2.2595 & 4.97E-3 & 2.1570 & 1.49E-2 & 1.6760  \\  
WENO-JS    & 40 & 2.83E-4 & 3.2361 & 8.51E-4 & 2.5439 & 3.75E-3 & 1.9854  \\
k = 3      & 80 & 1.53E-5 & 4.2131 & 5.48E-5 & 3.9567 & 3.20E-4 & 3.5528  \\ 
           & 160& 6.08E-7 & 4.6500 & 1.93E-6 & 4.8298 & 1.15E-5 & 4.7908  \\ 
           & 320& 1.80E-8 & 5.0792 & 5.56E-8 & 5.1168 & 3.35E-7 & 5.1071  \\
\hline  
           & 10 & 1.06E-2 &   --   & 1.87E-2 &    --  & 4.09E-2 &   --   \\  
           & 20 & 1.71E-3 & 2.6282 & 2.65E-3 & 2.8170 & 7.09E-3 & 2.5284 \\  
RBF-WENO-JS   & 40 & 2.53E-4 & 2.7592 & 6.44E-4 & 2.0428 & 2.62E-3 & 1.4336 \\
k = 3      & 80 & 1.40E-5 & 4.1748 & 4.82E-5 & 3.7417 & 2.77E-4 & 3.2434 \\ 
           & 160& 4.93E-7 & 4.8277 & 1.67E-6 & 4.8488 & 1.00E-5 & 4.7930 \\ 
           & 320& 1.46E-8 & 5.0755 & 4.83E-8 & 5.1130 & 2.88E-7 & 5.1170 \\ 
\hline 
\end{tabular}
\end{center} 
\label{table3}
\end{table}

\begin{figure}[h]
\begin{center}
\includegraphics[width=0.48\textwidth]{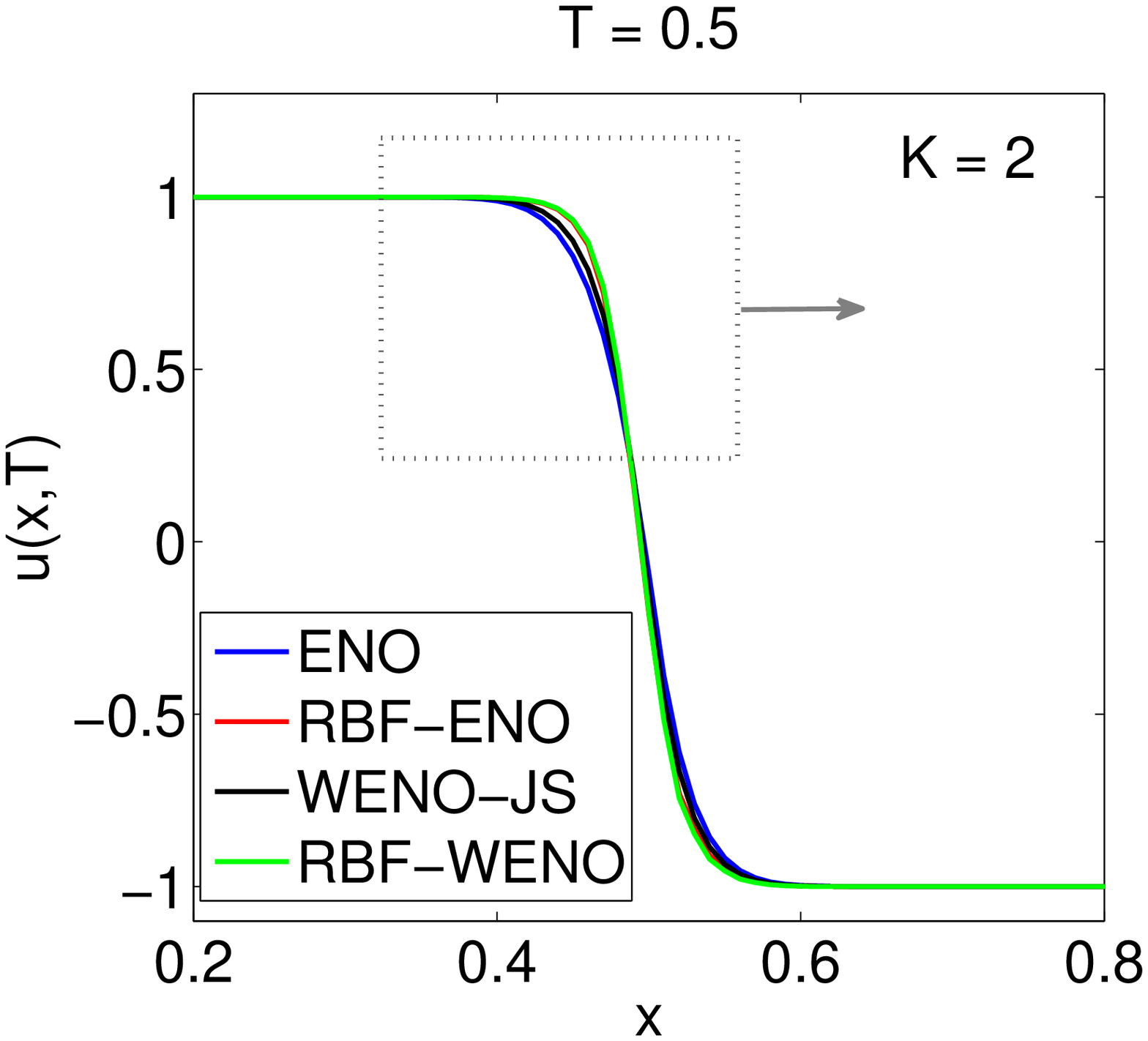}
\includegraphics[width=0.48\textwidth]{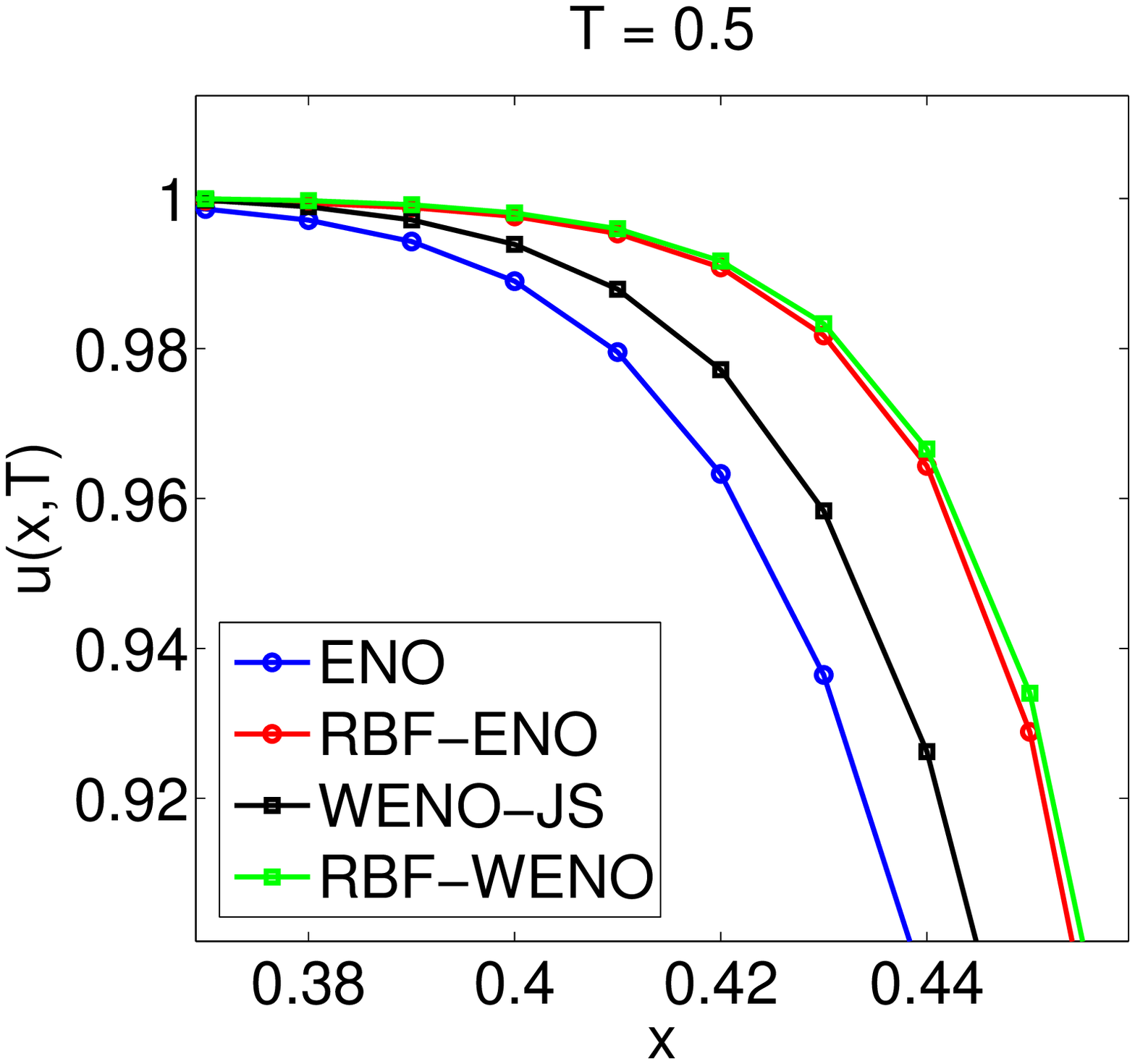}
\includegraphics[width=0.48\textwidth]{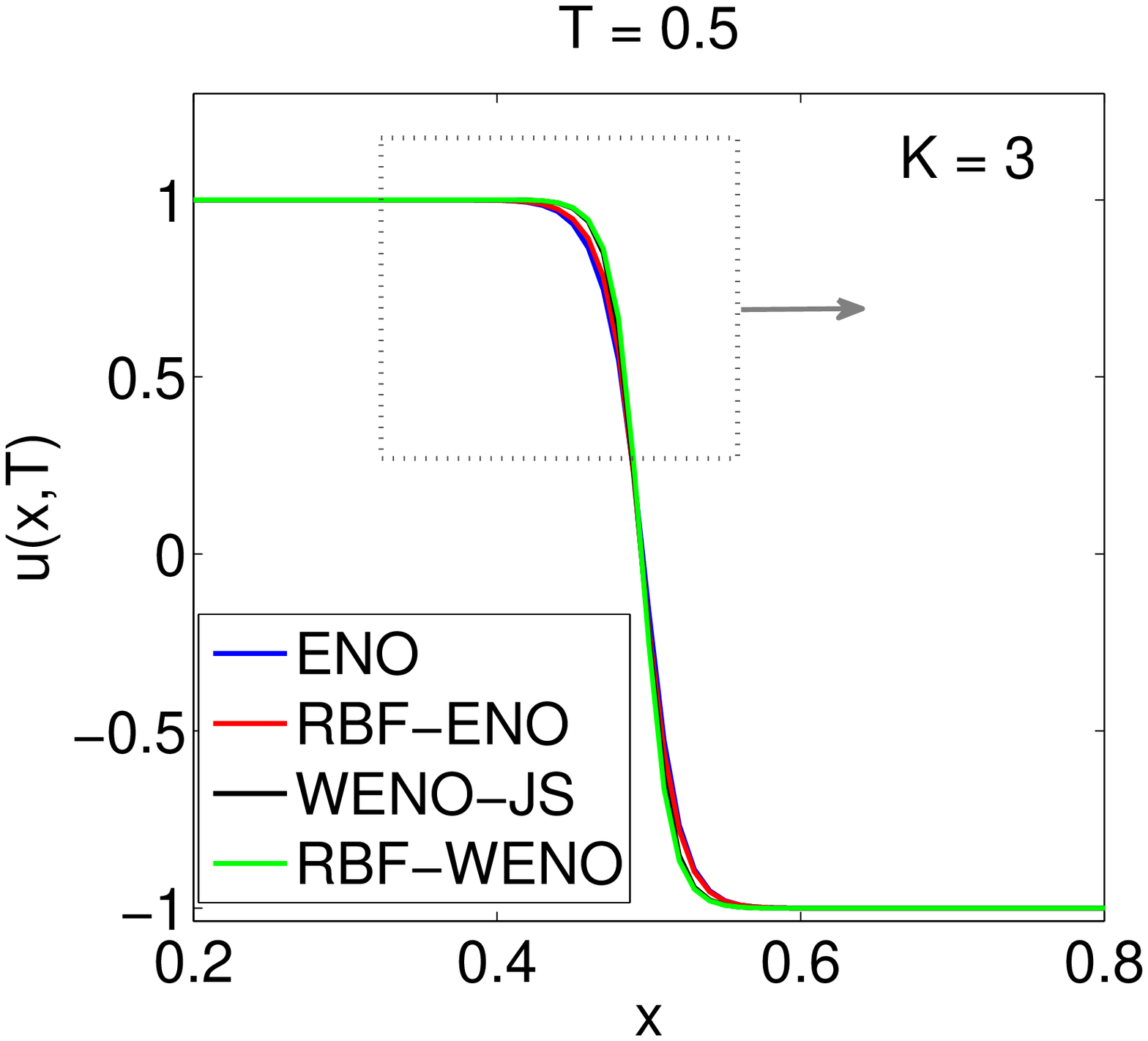}
\includegraphics[width=0.48\textwidth]{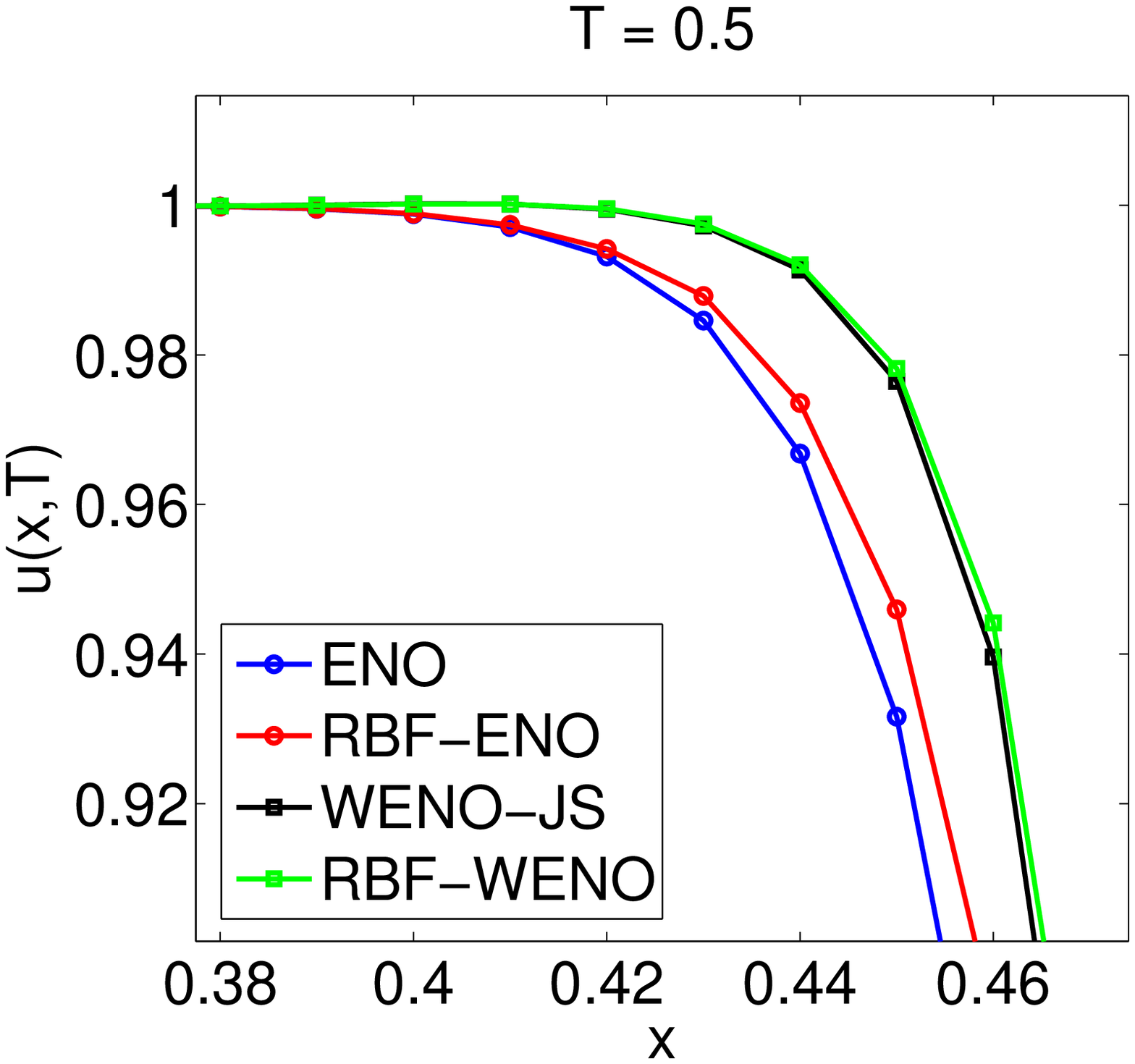}
\end{center}
\caption{Solutions to advection equation \eqref{example1} with discontinuous initial condition \eqref{example2}. $T = 0.5$. $k = 2$ (top) and $k = 3$ (bottom). $N = 200$.}
\label{figure2}

\bigskip
\begin{center}
\includegraphics[width=0.48\textwidth]{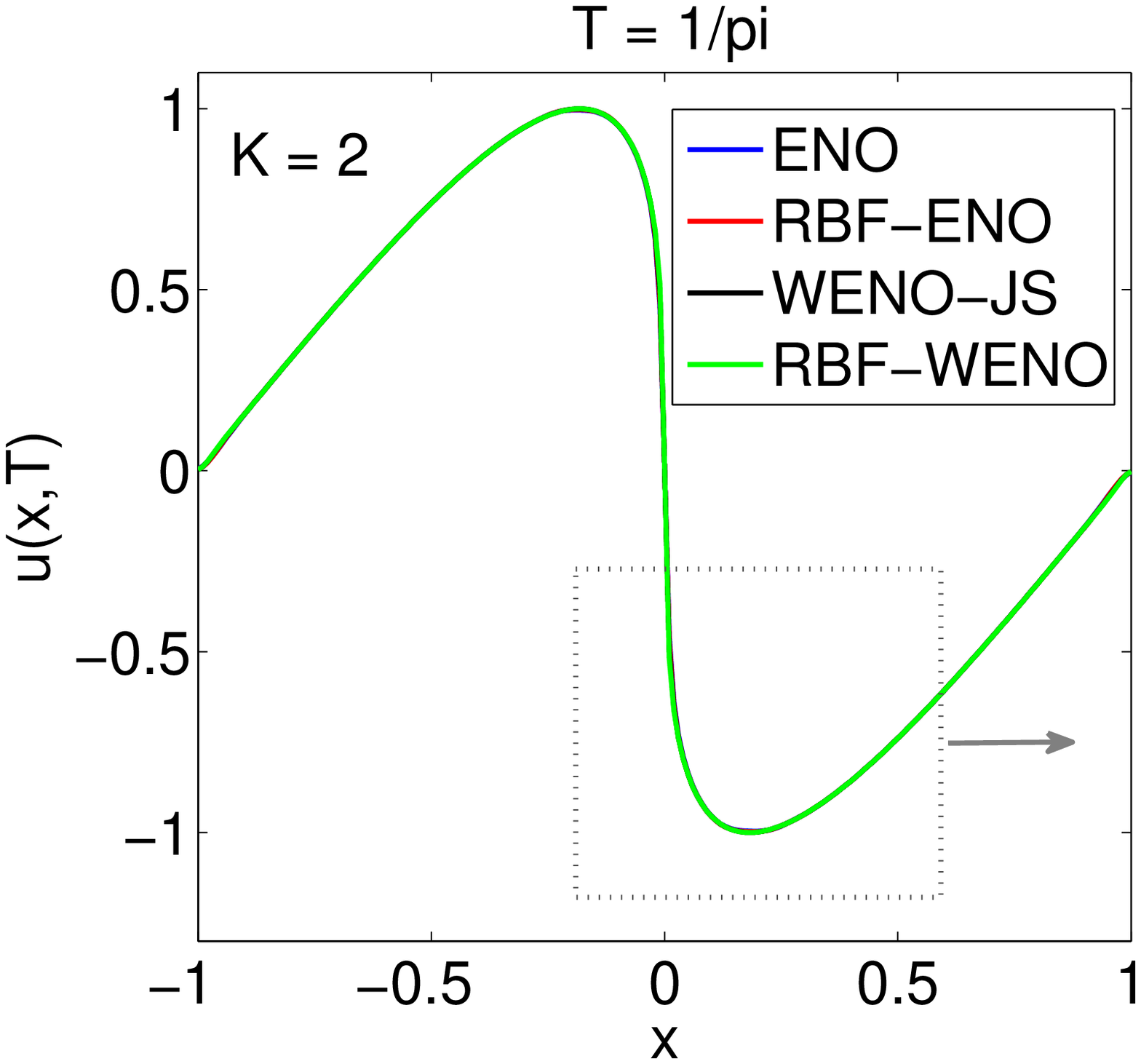}
\includegraphics[width=0.48\textwidth]{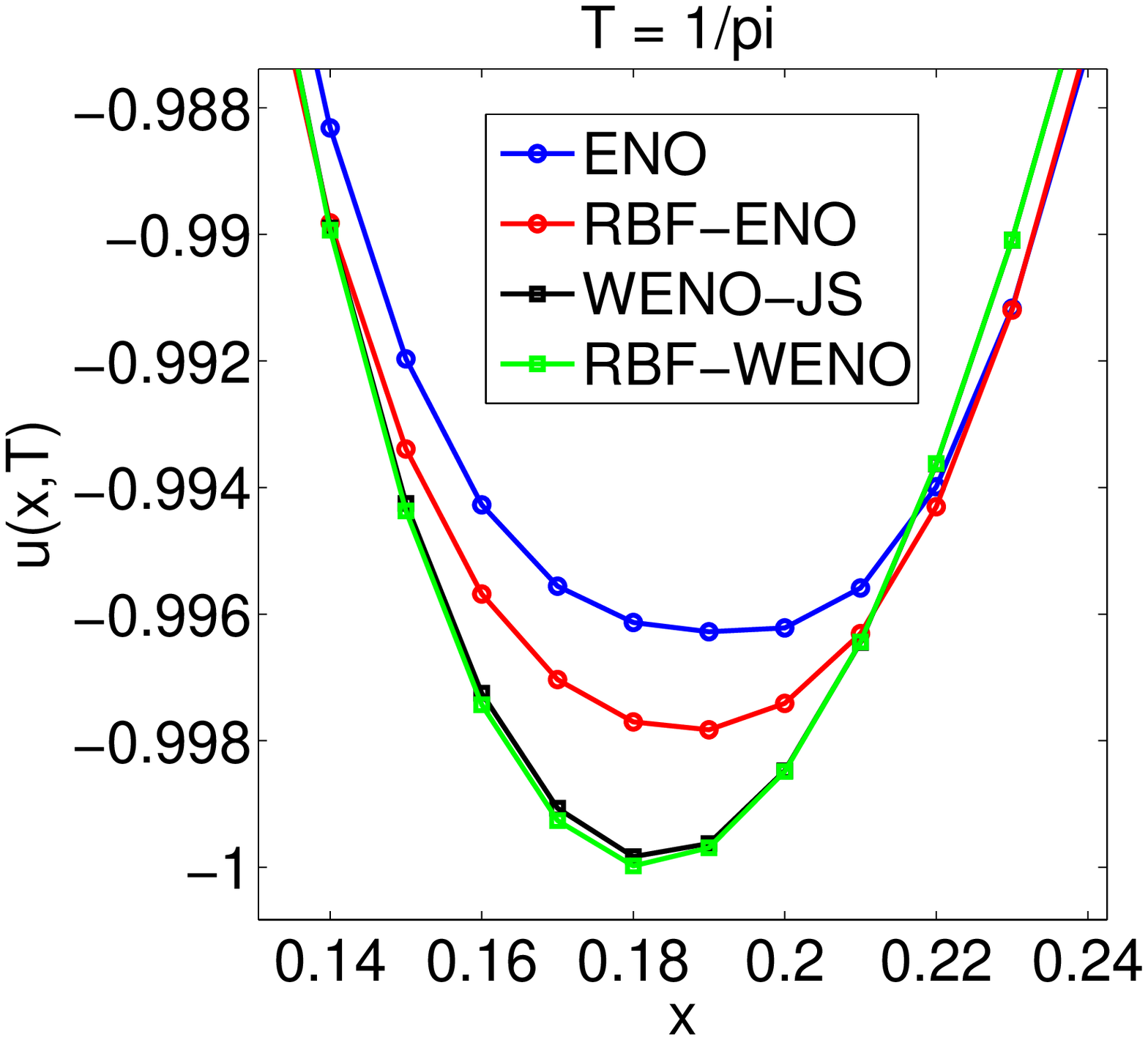}
\includegraphics[width=0.48\textwidth]{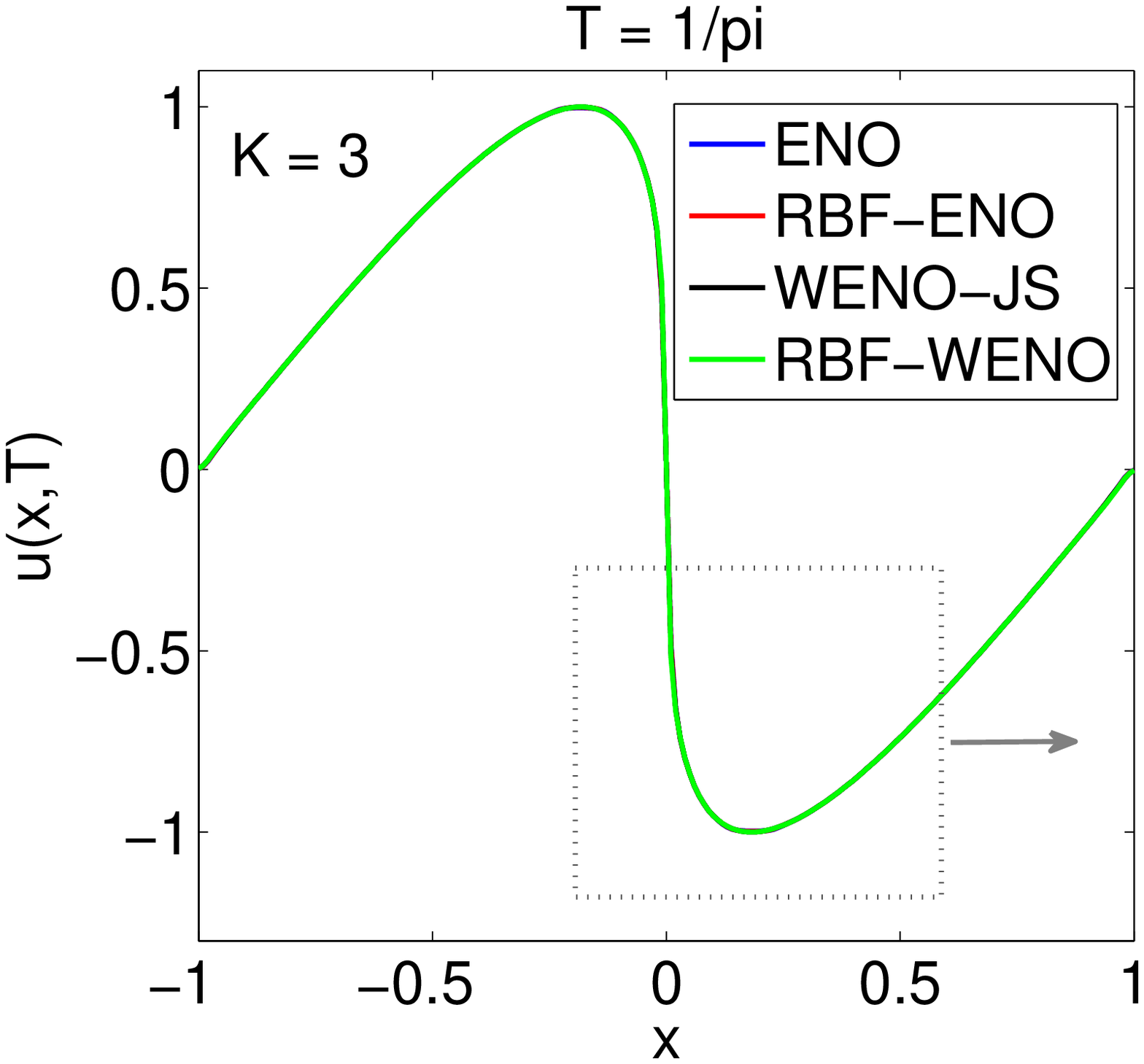}
\includegraphics[width=0.48\textwidth]{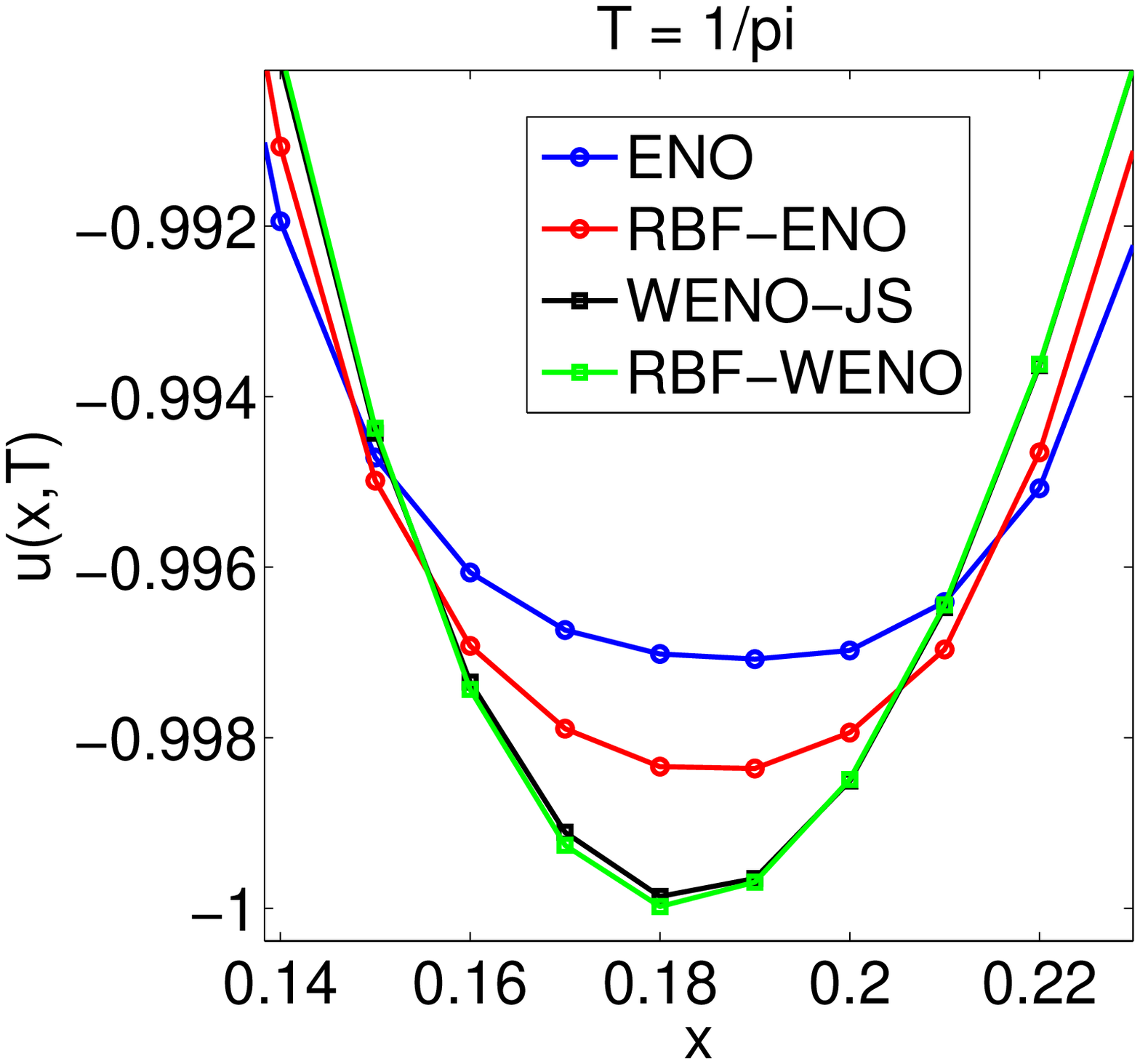}
\end{center}
\caption{(Color online). Solutions to Burger's equation \eqref{example3} with continuous initial condition \eqref{example3_example}. $T = \frac{1}{\pi}$. $k = 2$ (top) and $k = 3$ (bottom). $N = 200$.}
\label{figure3}
\end{figure}

Figure \ref{figure3} shows the solution profiles at the final time $T = {1\over \pi}$ when the shock forms with $k = 2$ and $k = 3$. The top figures show the solutions for $k = 2$ and the bottom for $k = 3$. For both cases, the RBF-WENO-JS and WENO-JS show similar results while the RBF-ENO solution is better than the regular ENO solution. 
%

\subsubsection{Euler equation I: Sod problem} 
\begin{figure}[h]
\begin{center}
\includegraphics[width=0.4\textwidth]{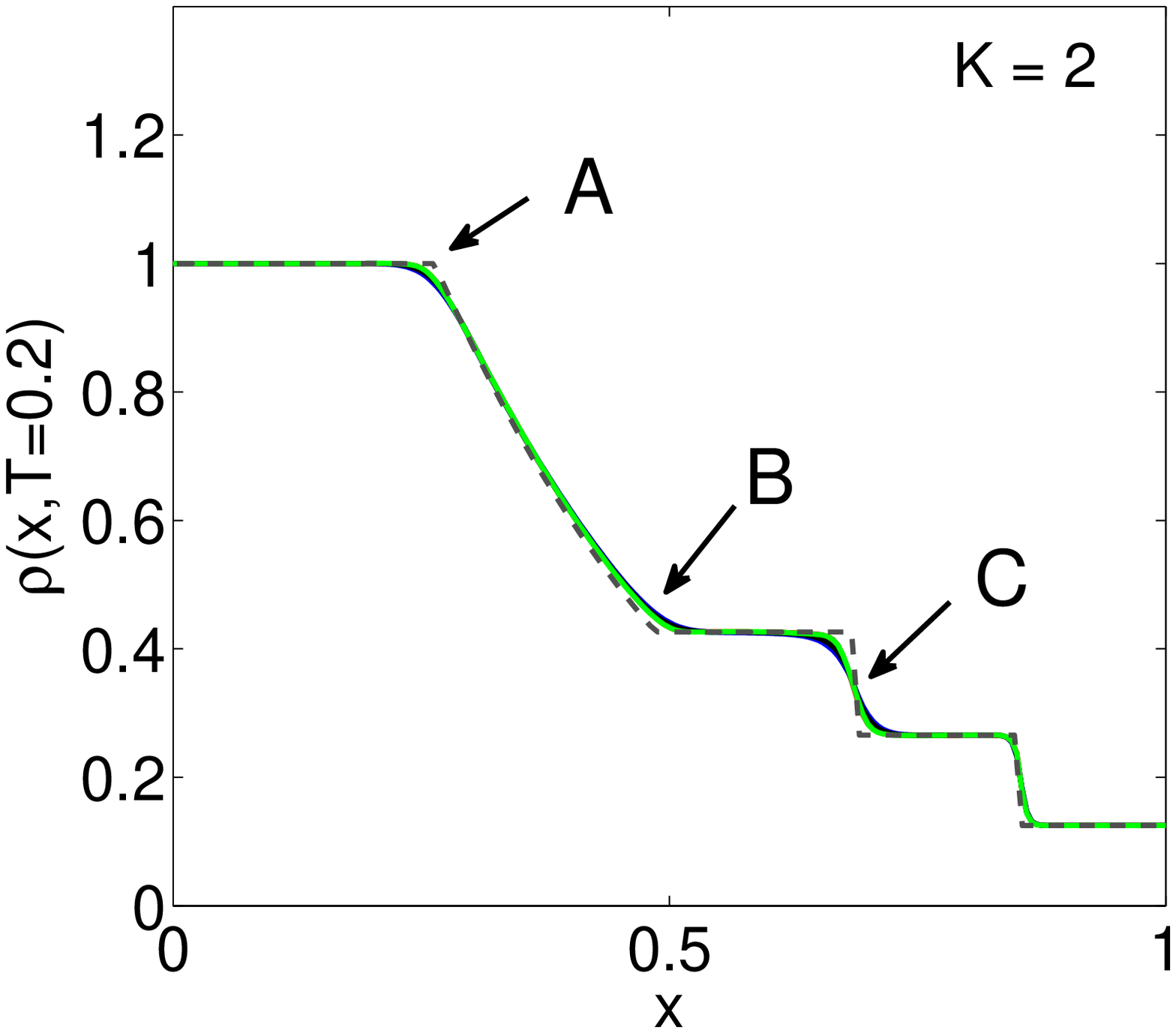}
\includegraphics[width=0.4\textwidth]{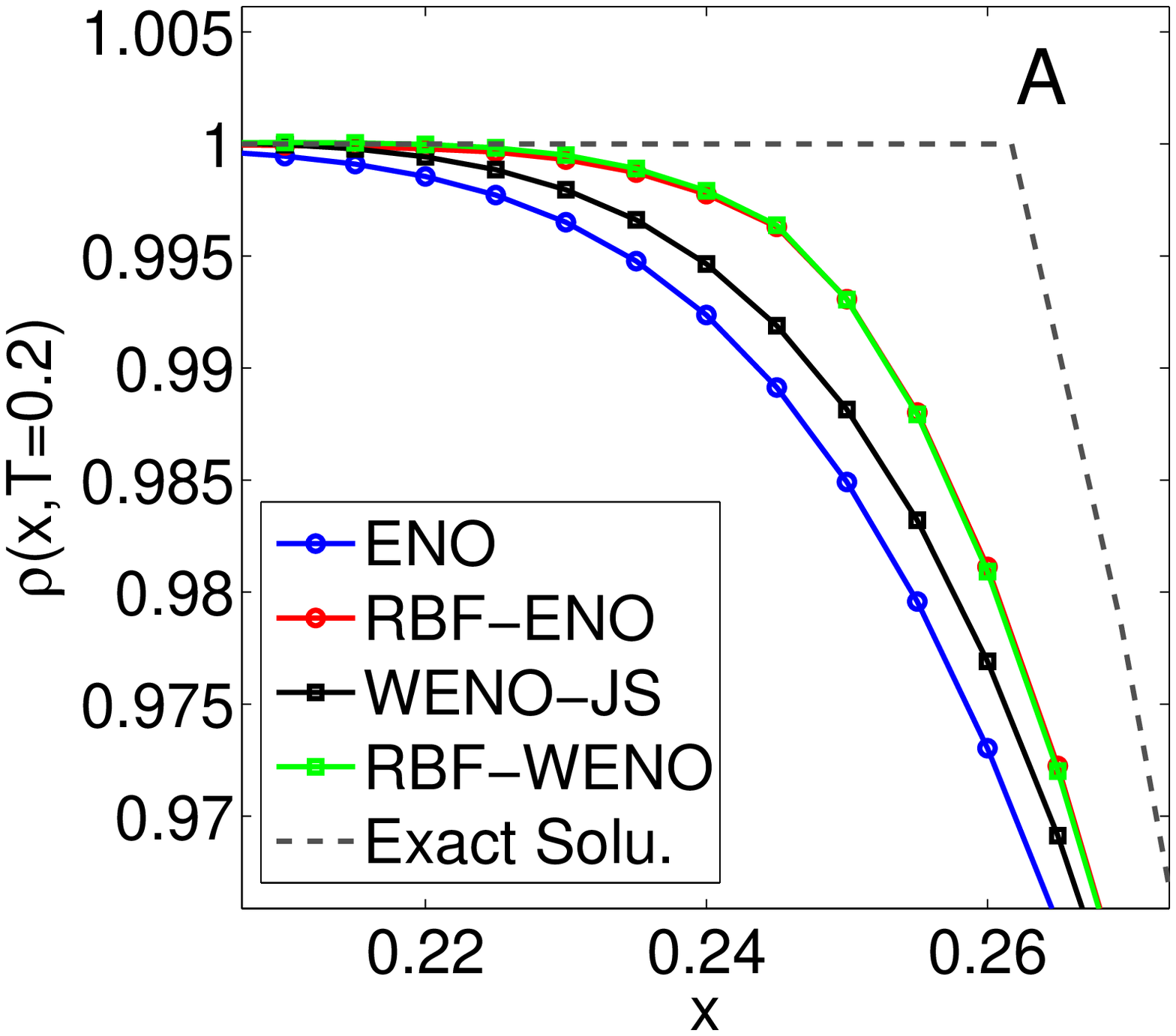}
\includegraphics[width=0.4\textwidth]{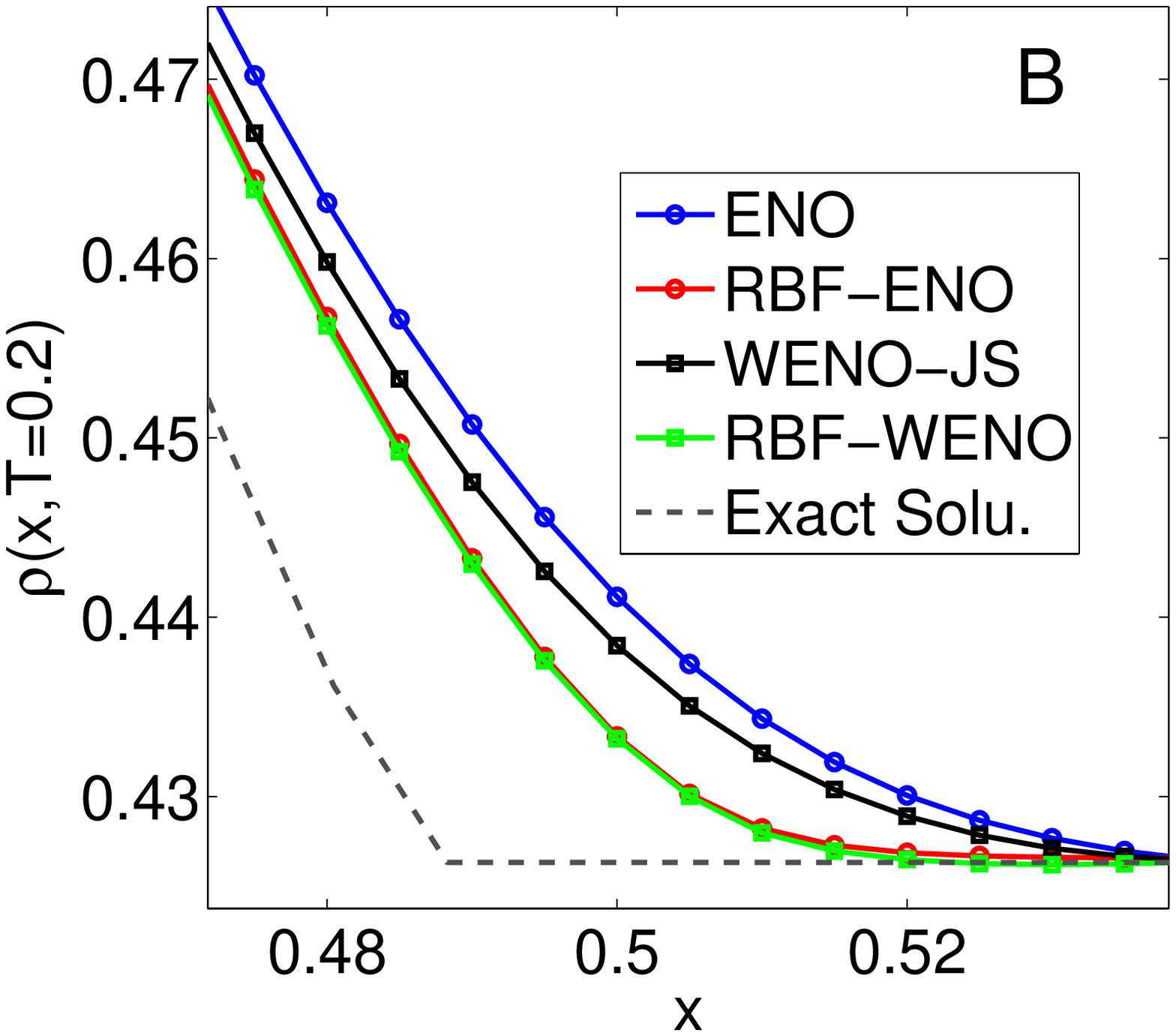}
\includegraphics[width=0.4\textwidth]{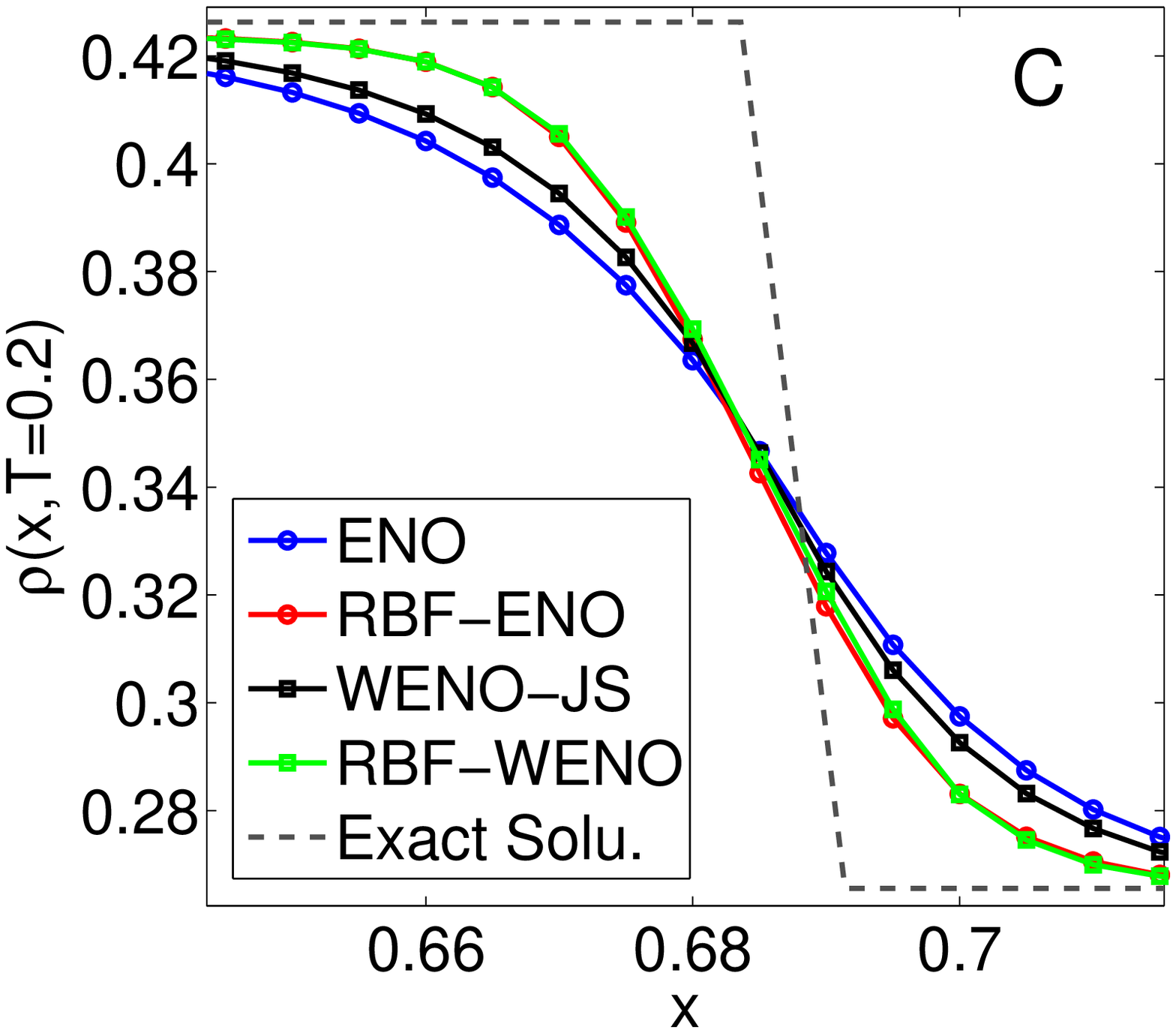}
\end{center}
\caption{Density profiles for Sod problem. $T = 0.2$. $N = 400$. $k = 2$. The dashed line is the exact solution. }
\label{sod}

\bigskip

\begin{center}
\includegraphics[width=0.4\textwidth]{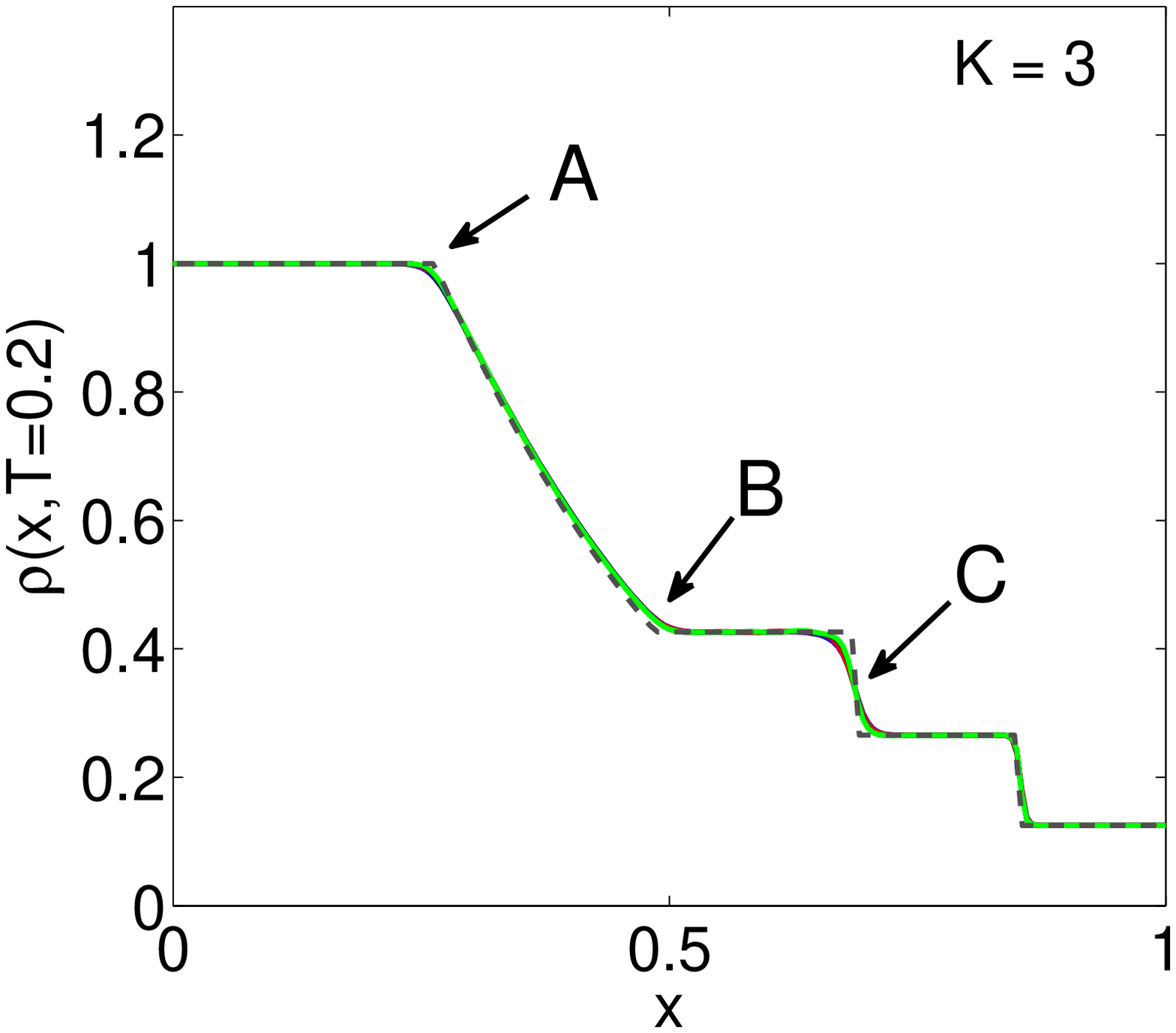}
\includegraphics[width=0.4\textwidth]{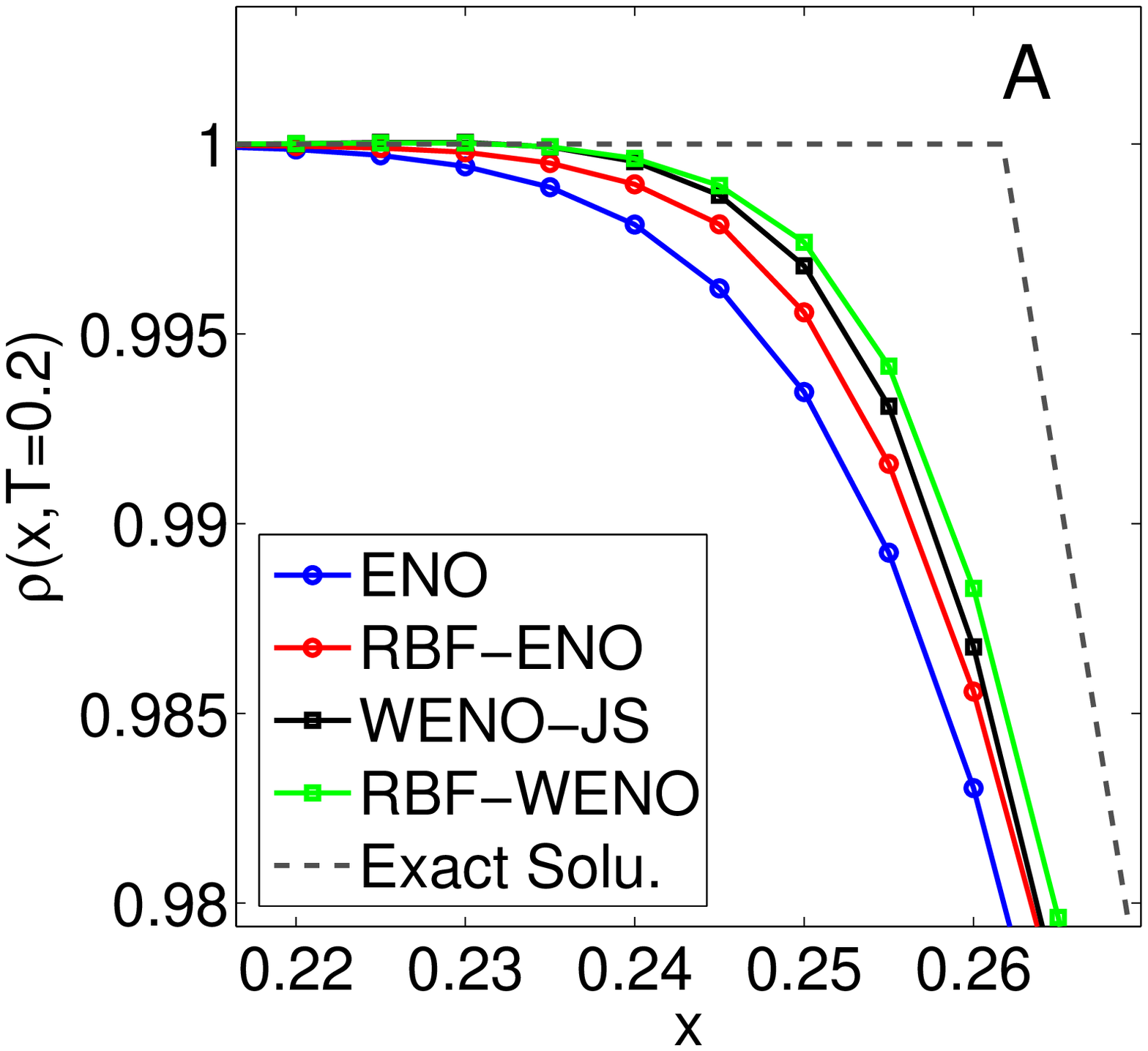}
\includegraphics[width=0.4\textwidth]{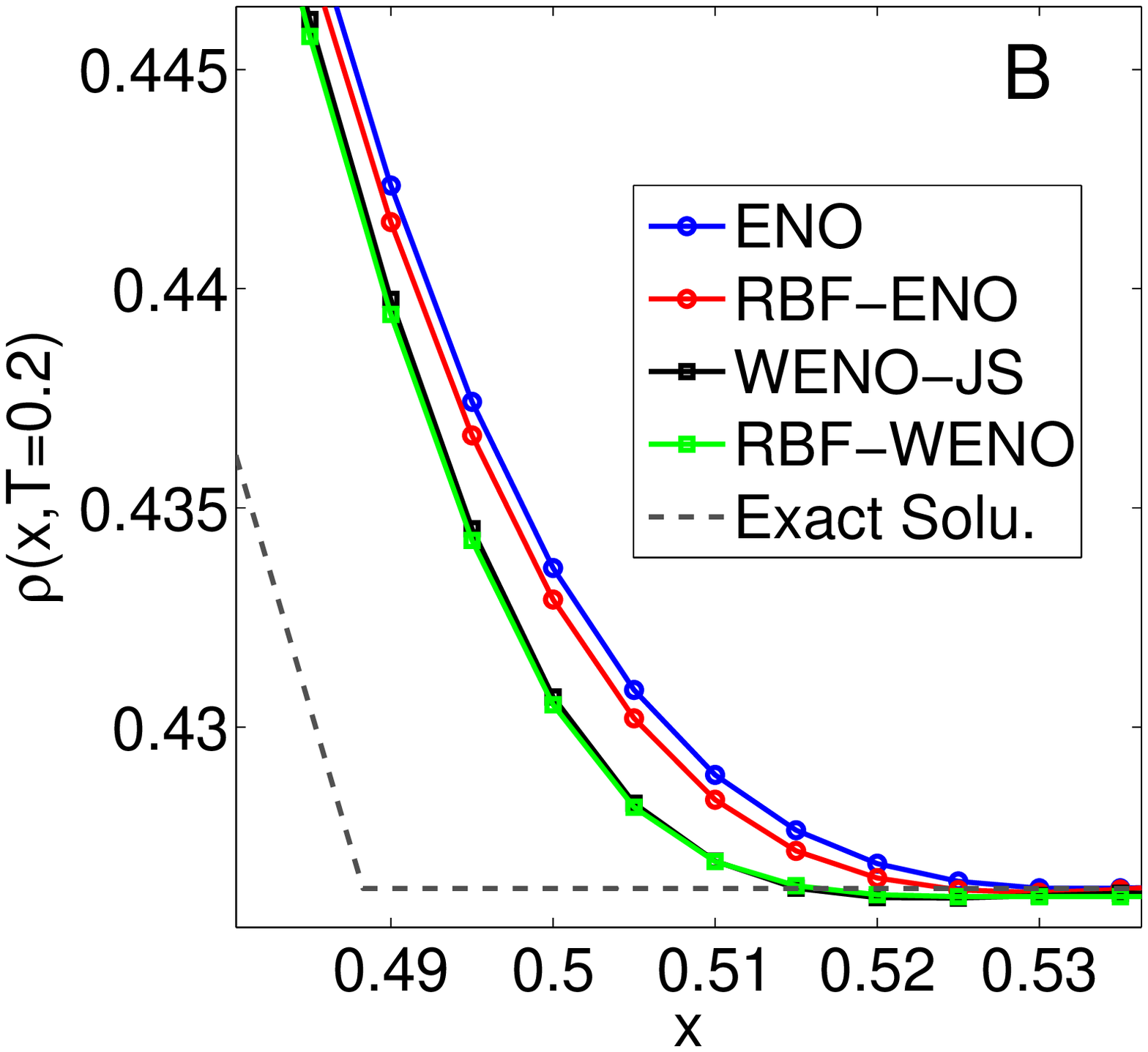}
\includegraphics[width=0.4\textwidth]{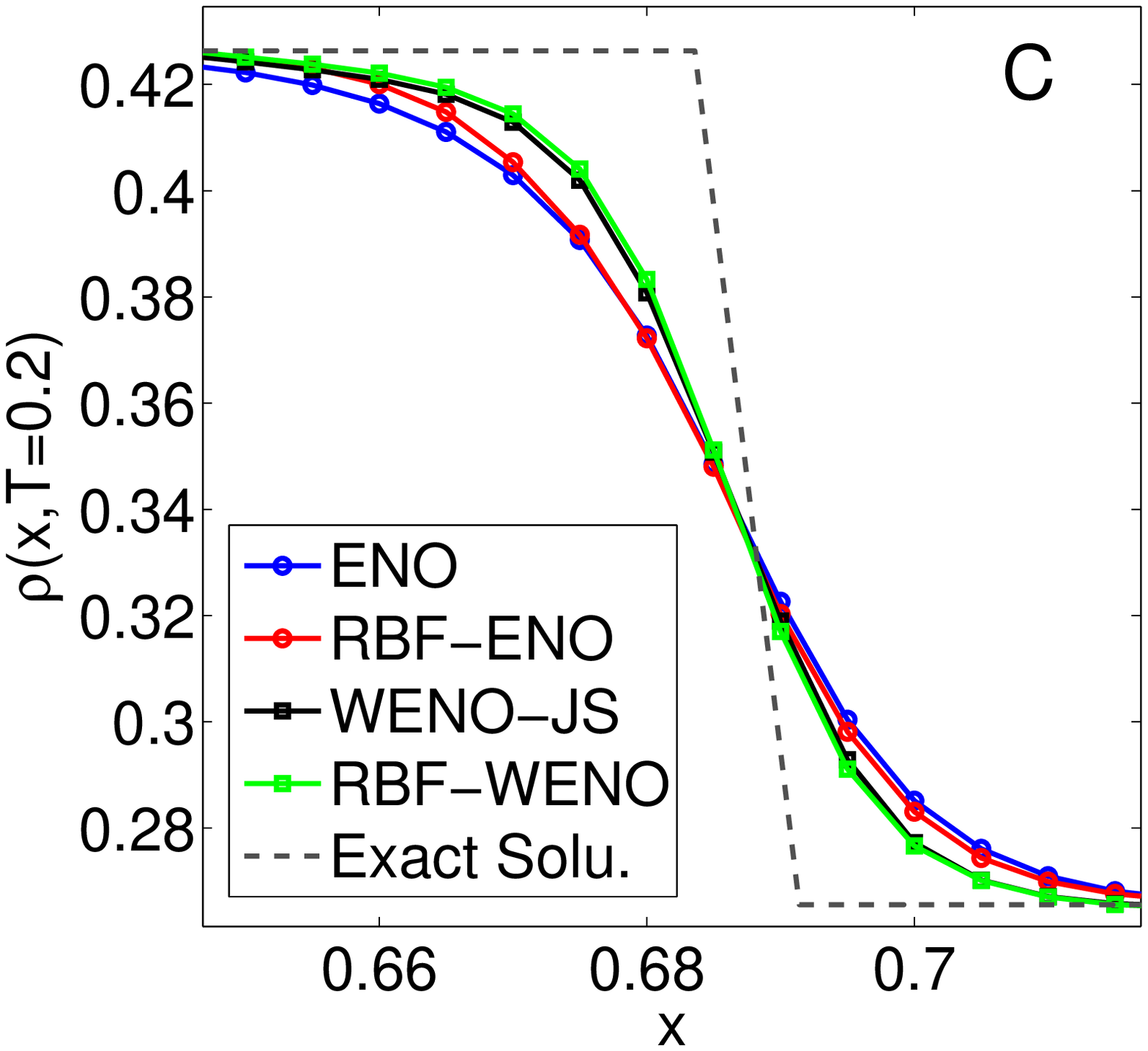}
\end{center}
\caption{Density profiles for Sod problem. $T = 0.2$. $N = 600$. $k = 3$. The dashed line is the exact solution. }
\label{sod2}
\end{figure}

For the system problem, we consider the one-dimensional Euler equations for gas dynamics
\begin{eqnarray}
    U_t + F(U)_x = 0, 
\label{euler}
\end{eqnarray}
where the conservative state vector $U$ and the flux function $F$ are given by 
$$
  U = (\rho, \rho u, E)^T, \quad F(U) = (\rho u, \rho u^2 + P, (E+P)u)^T.
$$
Here $\rho, u, P$ and $E$ denote density, velocity, pressure and total energy, respectively. The equation of state is given by 
$$
     P = (\gamma - 1) \left( E - {1\over 2} \rho u^2 \right), 
$$
where $\gamma = 1.4$ for the ideal gas. 
We consider the Sod shock tube problem with the initial condition 
$$
   (\rho, u, P ) = \left\{ \begin{array}{ll} (1, 0, 1) & x \le 0 \\ (0.125, 0, 0.1) & x > 0 \end{array} \right. . 
$$
The Sod problem is solved with $N = 200$ and the CFL number $C = 0.1$. Figures \ref{sod} and \ref{sod2} show solutions by each method at the final time $T = 0.2$ for $k = 2$ and $k = 3$, respectively. The top left figure shows the global solution in the entire domain while the rest shows the solution profiles in the areas specified in the top left figure. As in the previous examples, for $k = 2$ the RBF-ENO and RBF-WENO-JS methods are better performed than the regular ENO and WENO-JS methods. The solution profiles are much sharper near the non-smooth area.  For $k = 3$, the RBF-WENO-JS solution is similar to the WENO-JS solution while the RBF-ENO solution is sharper than the regular ENO solution.  
%
%

\subsubsection{Euler equation II: Lax problem} 
\begin{figure}[h]
\begin{center}
\includegraphics[width=0.4\textwidth]{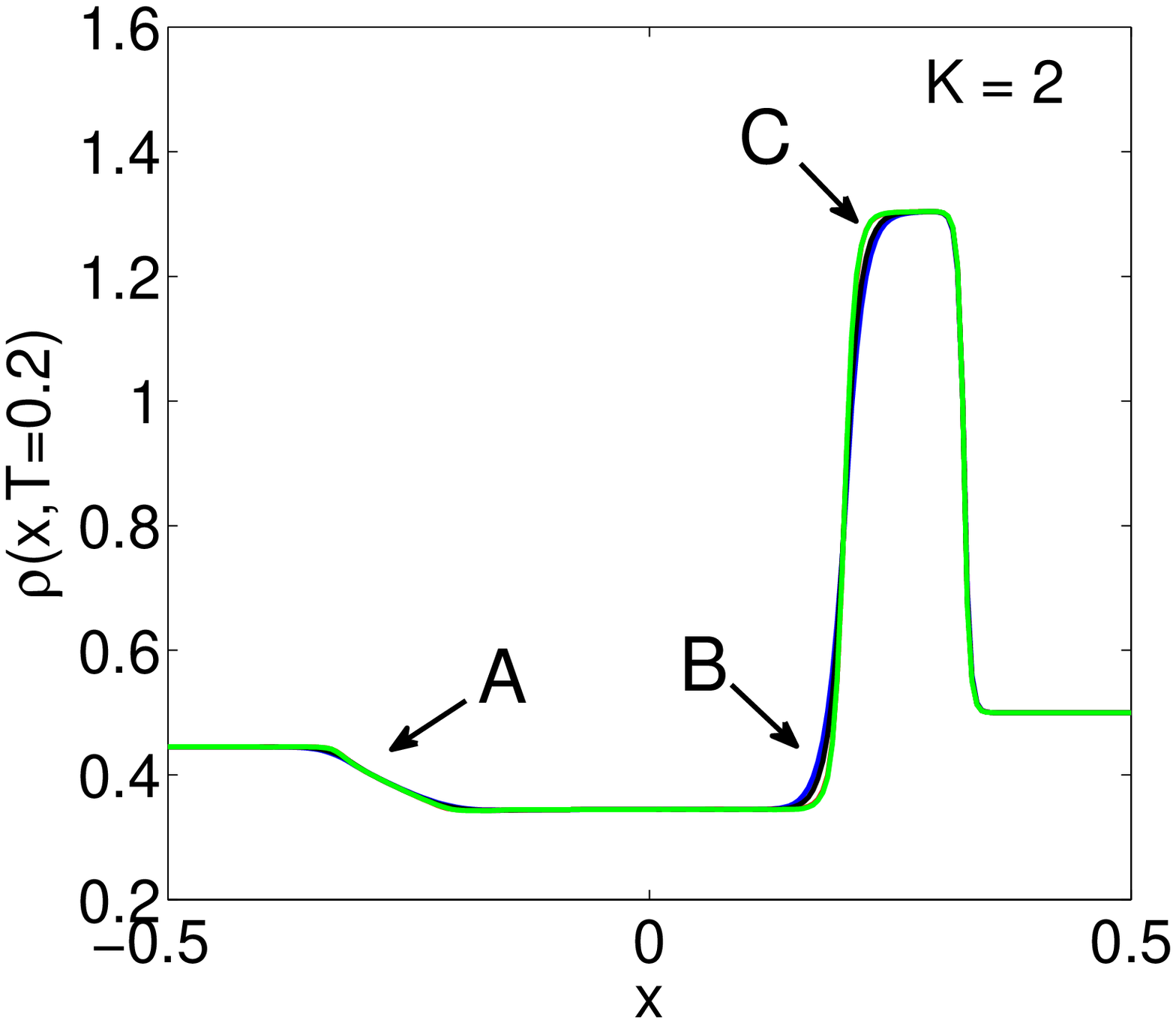}
\includegraphics[width=0.4\textwidth]{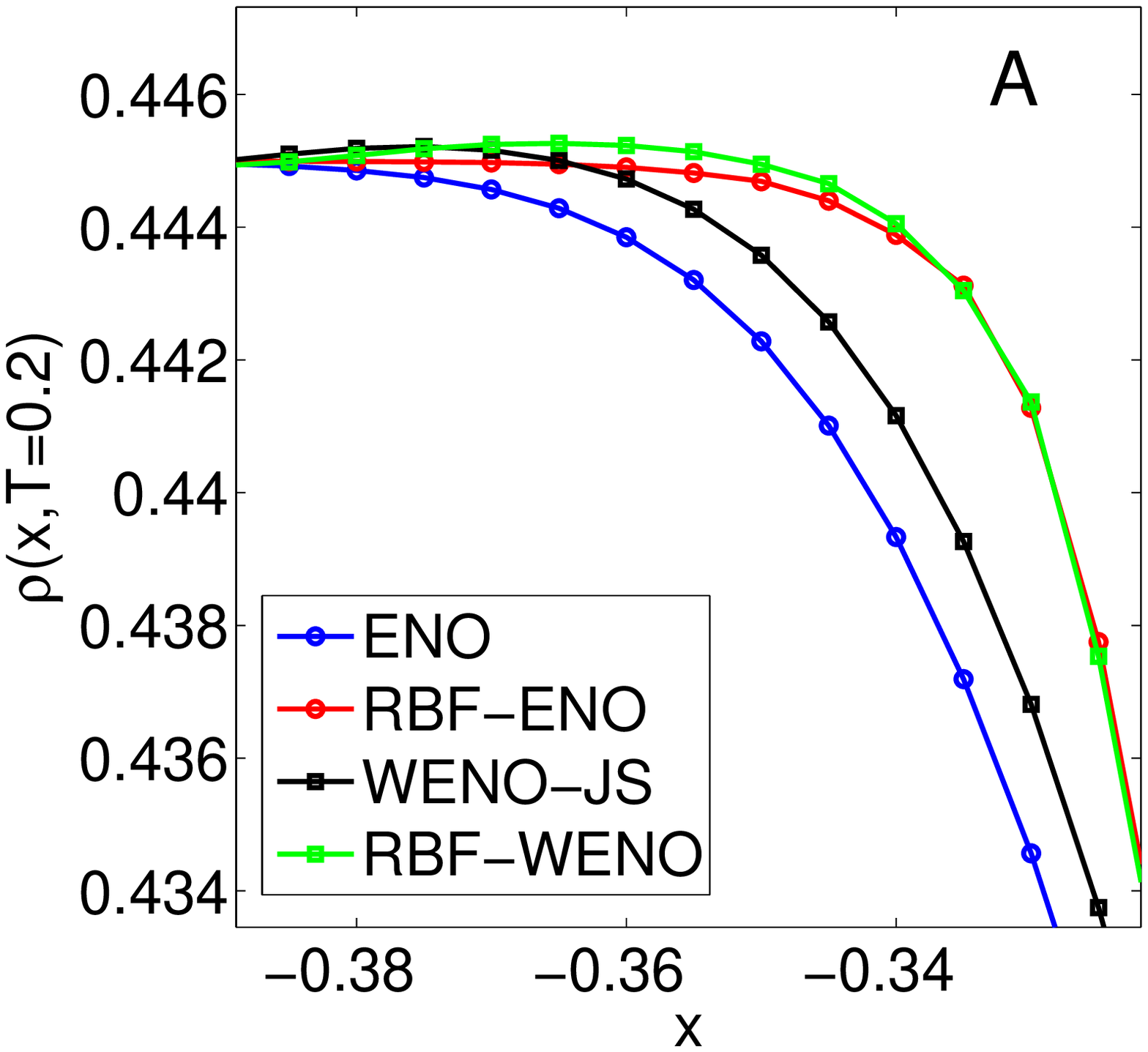}
\includegraphics[width=0.4\textwidth]{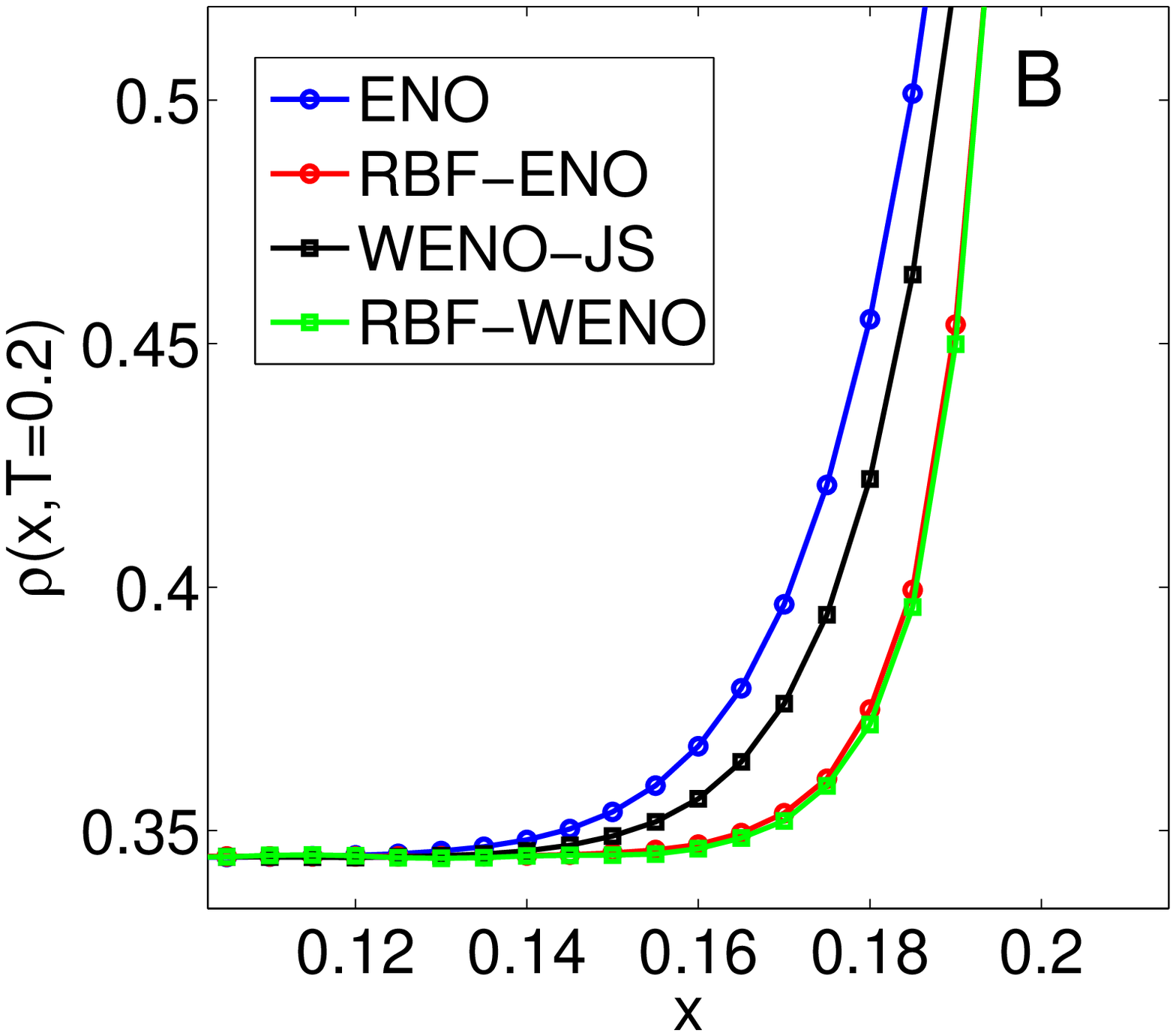}
\includegraphics[width=0.4\textwidth]{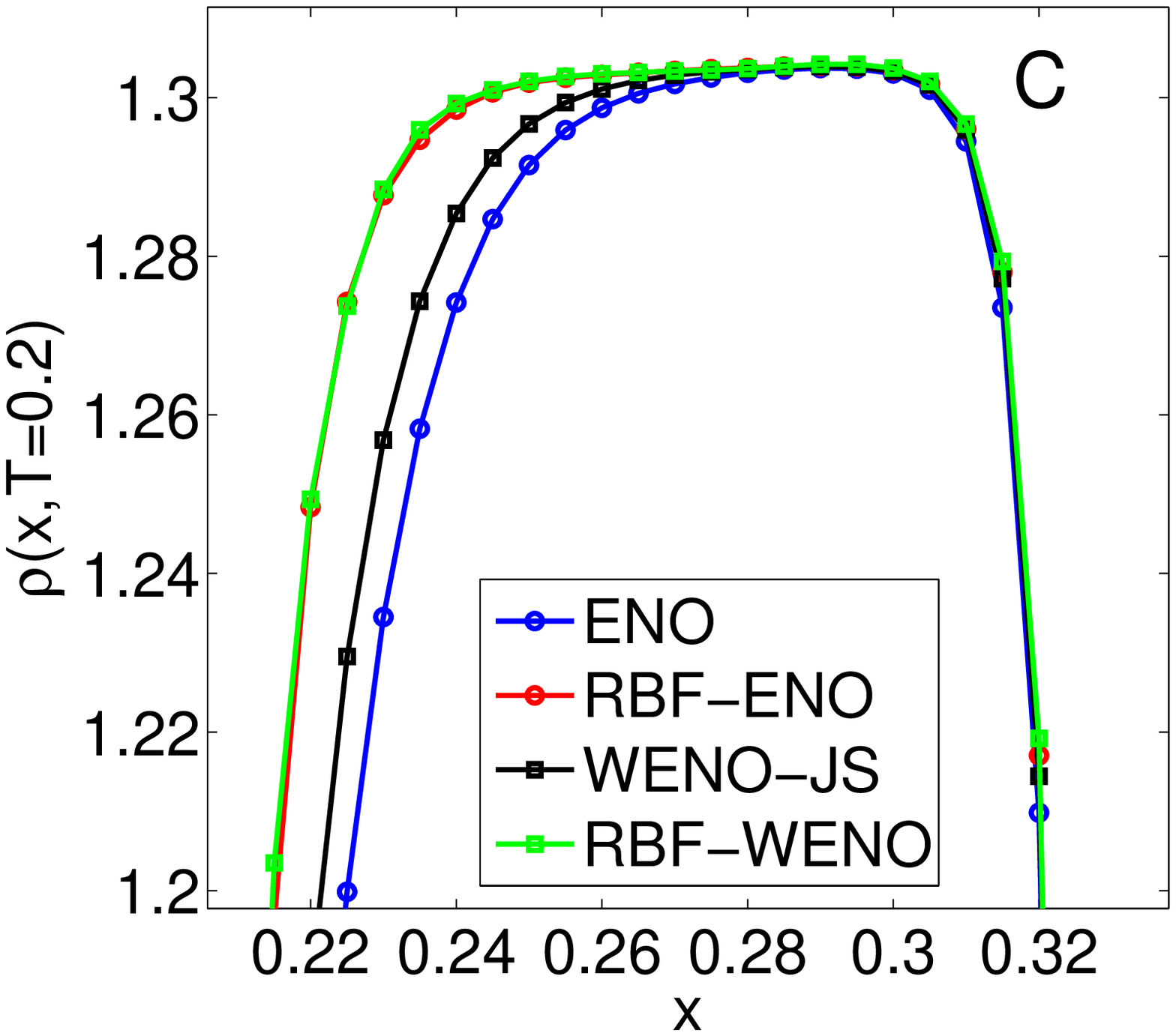}
\end{center}
\caption{Density profiles for Lax problem. $T = 0.13$. $N = 200$. $k = 2$.}
\label{lax1}

\bigskip

\begin{center}
\includegraphics[width=0.4\textwidth]{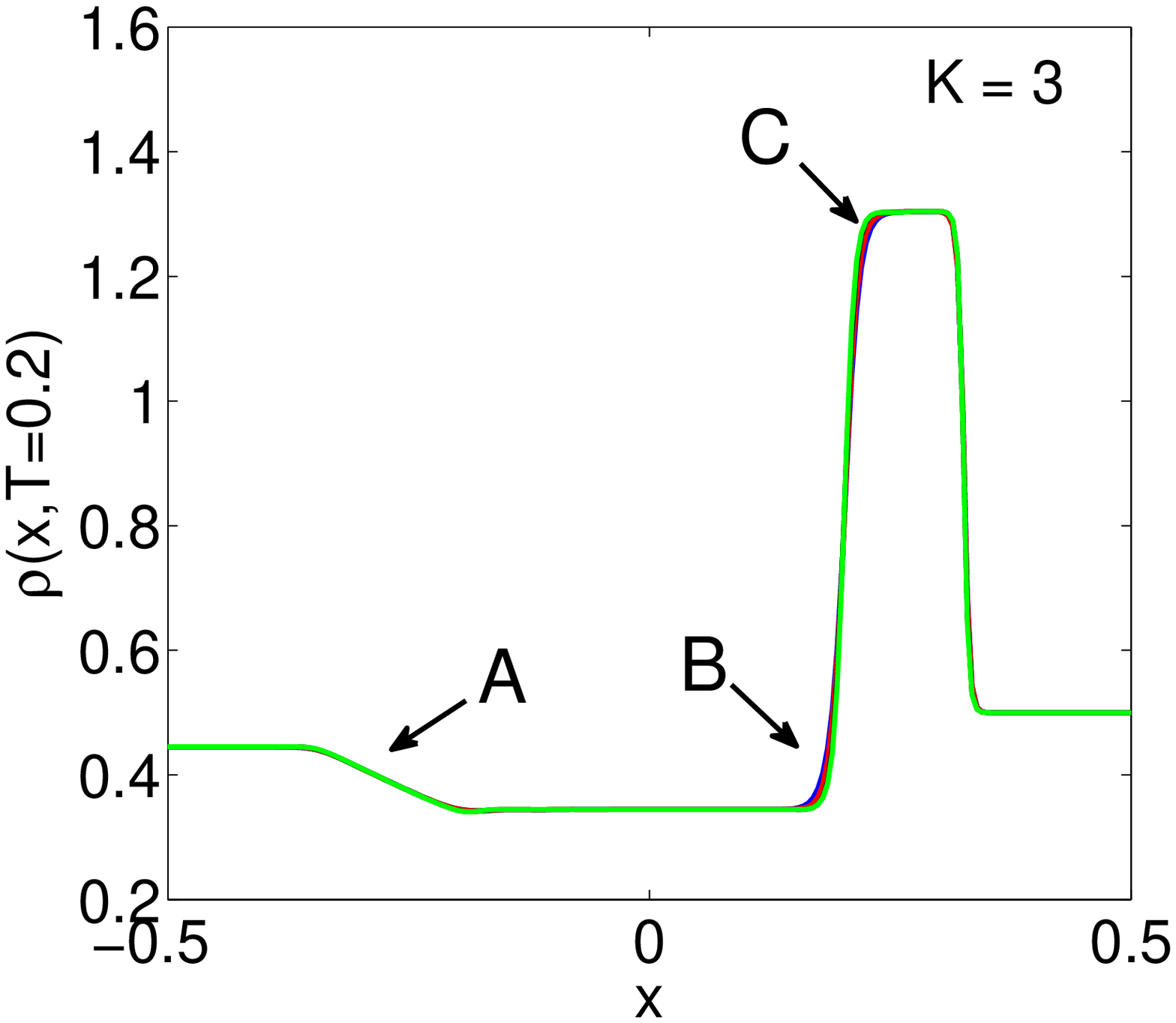}
\includegraphics[width=0.4\textwidth]{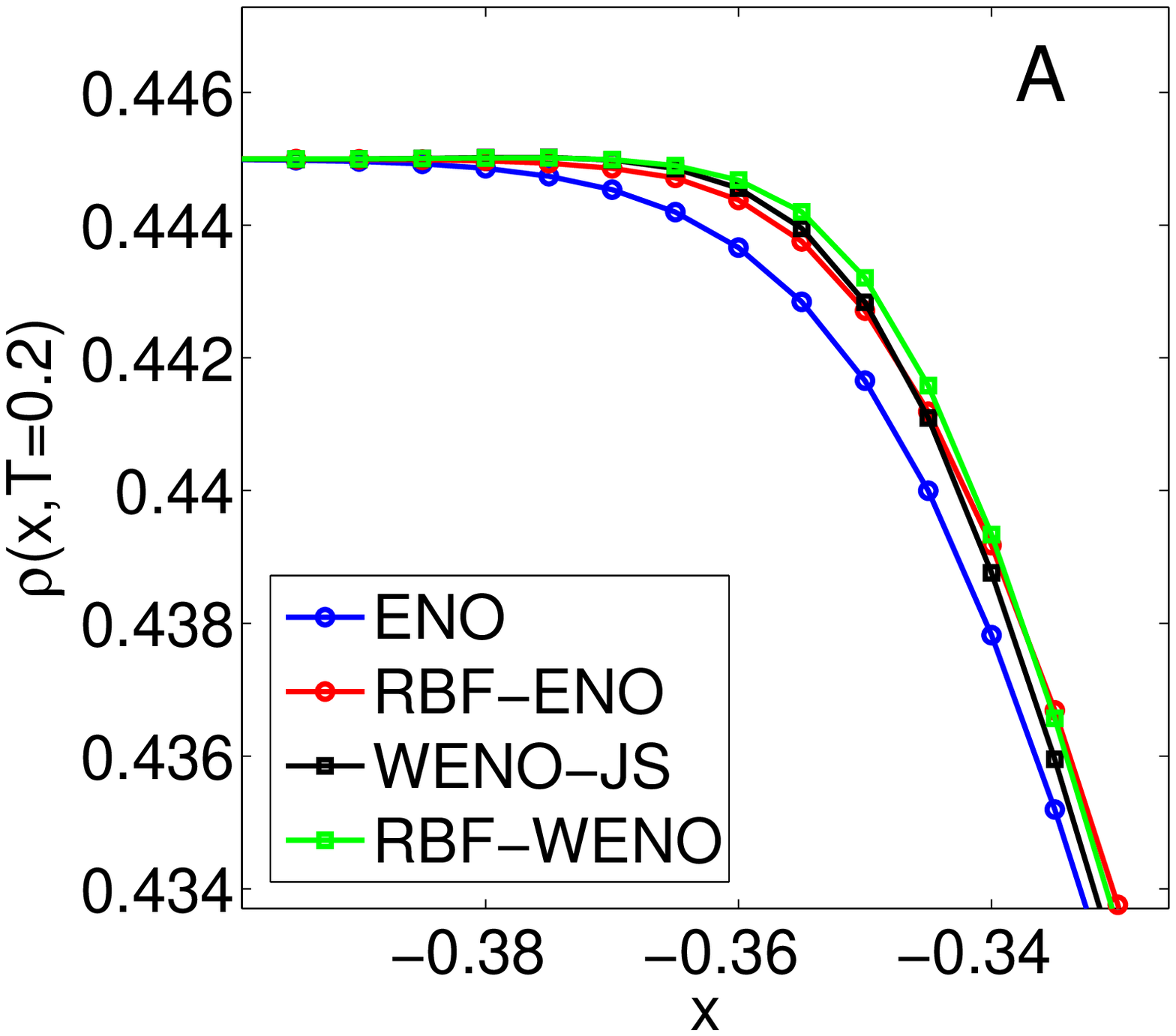}
\includegraphics[width=0.4\textwidth]{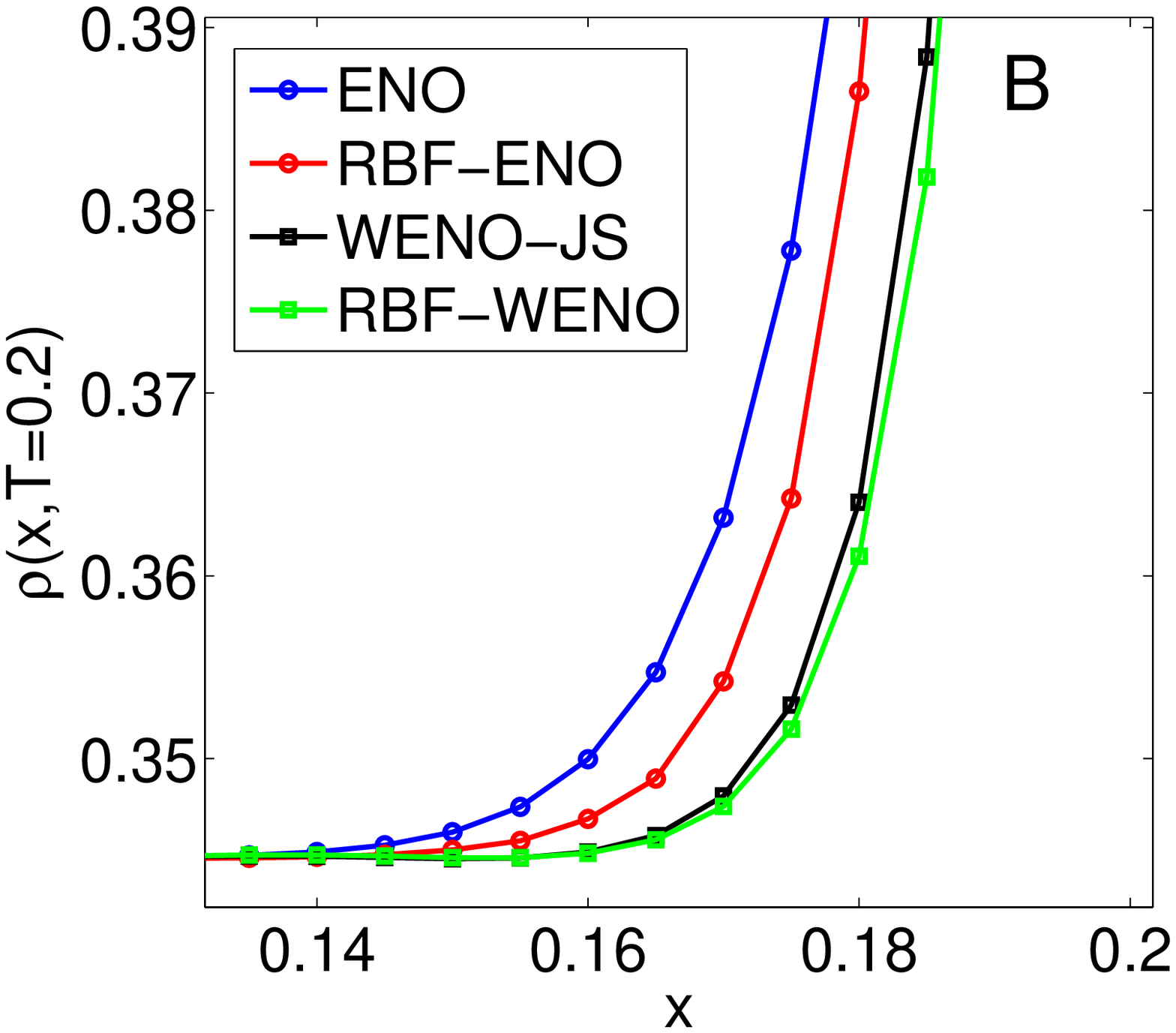}
\includegraphics[width=0.4\textwidth]{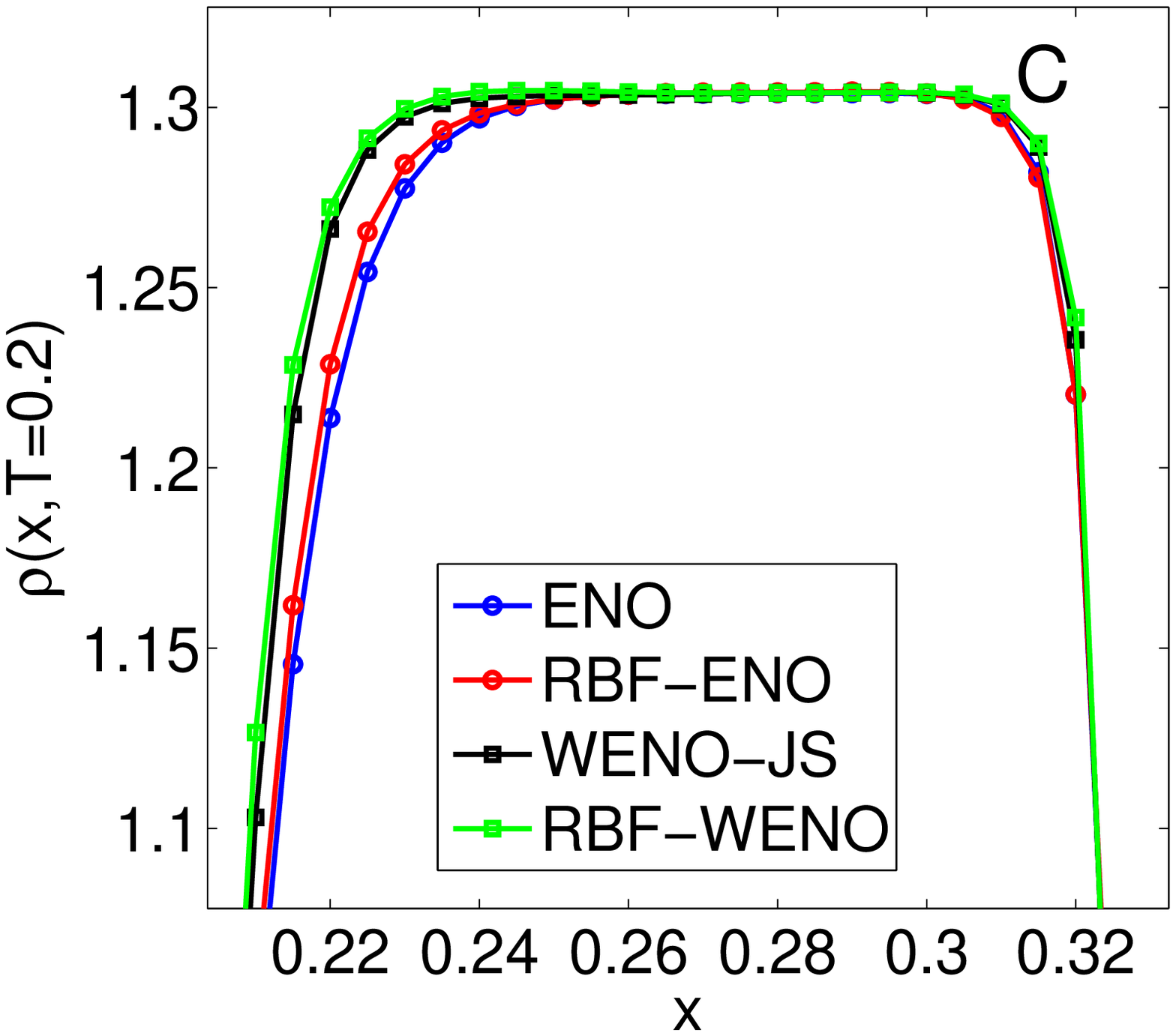}
\end{center}
\caption{Density profiles for Lax problem. $T = 0.13$. $N = 200$. $k = 3$.}
\label{lax2}
\end{figure}

The Lax problem has the following initial conditions:
$$
   (\rho, u, P ) = \left\{ \begin{array}{ll} (0.445, 0.698, 3.528) & x \le 0 \\ (0.5, 0, 0.571) & x > 0 \end{array} \right. 
$$
with the Dirichlet boundary conditions at $x=0.5$ and $x=-0.5$. The solution was integrated until $t=0.13$. Figures \ref{lax1} and \ref{lax2} show the solution profiles at the final time $T = 0.13$ with $N = 200$ for $k = 2$ and $k = 3$, respectively. The top left figure shows the global solution and the rest the local solution profile in the areas specified in the top left figure. For $k = 2$, as in the previous example, the RBF-ENO and RBF-WENO-JS solutions are better than the regular ENO and WENO solutions near the non-smooth area. For $k = 3$, the RBF-WENO-JS solution is slightly better than the WENO solution and the RBF-ENO solution is better than the regular ENO solution. Since the RBF-WENO-JS and WENO-JS methods both yield about $5$th order accuracy while the RBF-ENO solution is $3$rd order or higher but less than $5$th order, it is expected that the RBF-WENO-JS or WENO-JS solutions are better than the ENO and RBF-ENO solutions. 


\subsection{2D Numerical examples}
Now we consider the two-dimensional hyperbolic conservation law below
$$ 
u_t(x,y,t)+f_x(u(x,y,,t))+g_y(u(x,y,,t)) = 0, 
$$
with the appropriate initial and boundary conditions where $f$ and $g$ are the flux functions. As in 1D case, the ODE system for the discretized version in the cell $I_{ij} = [x_{i-{1\over 2}}, x_{i + {1\over 2}}] \times [y_{i-{1\over 2}}, y_{i + {1\over 2}}]$ is given by 
\begin{eqnarray}
  \frac{du_{ij}(t)}{dt} =  -\frac{1}{\Delta x} \left[ \hat{f}_{i+{1\over 2},j} - \hat{f}_{i-{1\over 2},j} \right] -\frac{1}{\Delta y} \left[ \hat{g}_{i,j+{1\over 2}} - \hat{h}_{i,j-{1\over 2}} \right],
\end{eqnarray} 
where $u_{ij}(t)$ is the numerical approximation to $u(x_i,y_j,t)$ and ${\hat f}$ and $\hat g$ are the numerical fluxes for the flux functions $f$ and $g$, respectively. Notice that unlike the finite volume method, the finite difference method does not involve any quadrature points for the numerical integration of the flux functions. With the finite volume non-polynomial ENO/WENO methods for 2D problems, the accuracy of the local reconstruction can be enhanced at certain quadrature points. But since there is no quadrature point involved, the finite difference RBF-WENO-JS method maintains the same order of convergence although the accuracy of the RBF-WENO-JS method is enhanced. 

For the 2D numerical experiments, we consider the double Mach reflection problem described in \cite{Woodward} with the domain $(x,y) \in [0, 4] \times [0, 1]$.  The simulation is run until the final time $T=0.2$. At $T = 0.2$,  shocks form near the right edge shown as an almost straight line in the density contour in the following figures. 
Figures \ref{DMR1} and \ref{DMR11} show the density contours with various methods. The left column of figures shows the density contours with the regular ENO/WENO methods and the right column with the RBF-ENO/WENO-JS methods. These contours look similar with slight differences. For a detailed comparison, Figs. \ref{DMR2} and \ref{DMR22} show a slice image of Figs. \ref{DMR1} and \ref{DMR11} at $y=0.5$. We zoom in different regions for the case of $N=160$ and $M=40$ (Fig. \ref{DMR2}) and $N=320$ and $M=80$ (Fig. \ref{DMR22}) to compare all the methods. Here $N$ and $M$ denote the total number of grids in $x$ and $y$ directions, respectively. With $k=2$, we can see that the RBF-ENO solutions are comparable or even better than the regular WENO solutions. This is because when the RBF-ENO and the regular WENO solutions are of the same order, the RBF-ENO solution is sharper and less dissipative. The RBF-WENO-JS solution is also slightly better and sharper than the regular WENO-JS solution. With $k=3$, the regular WENO solution is better than the RBF-ENO solution  since the accuracy of the RBF-ENO method is lower than the regular WENO method. However, the RBF-ENO solution is still sharper than the regular ENO solution and has a similar shape to the WENO solution. We also observe that the RBF-WENO-JS method performs best among all those methods.  

\begin{figure}[h]
\begin{center}
\includegraphics[width=0.8\textwidth]{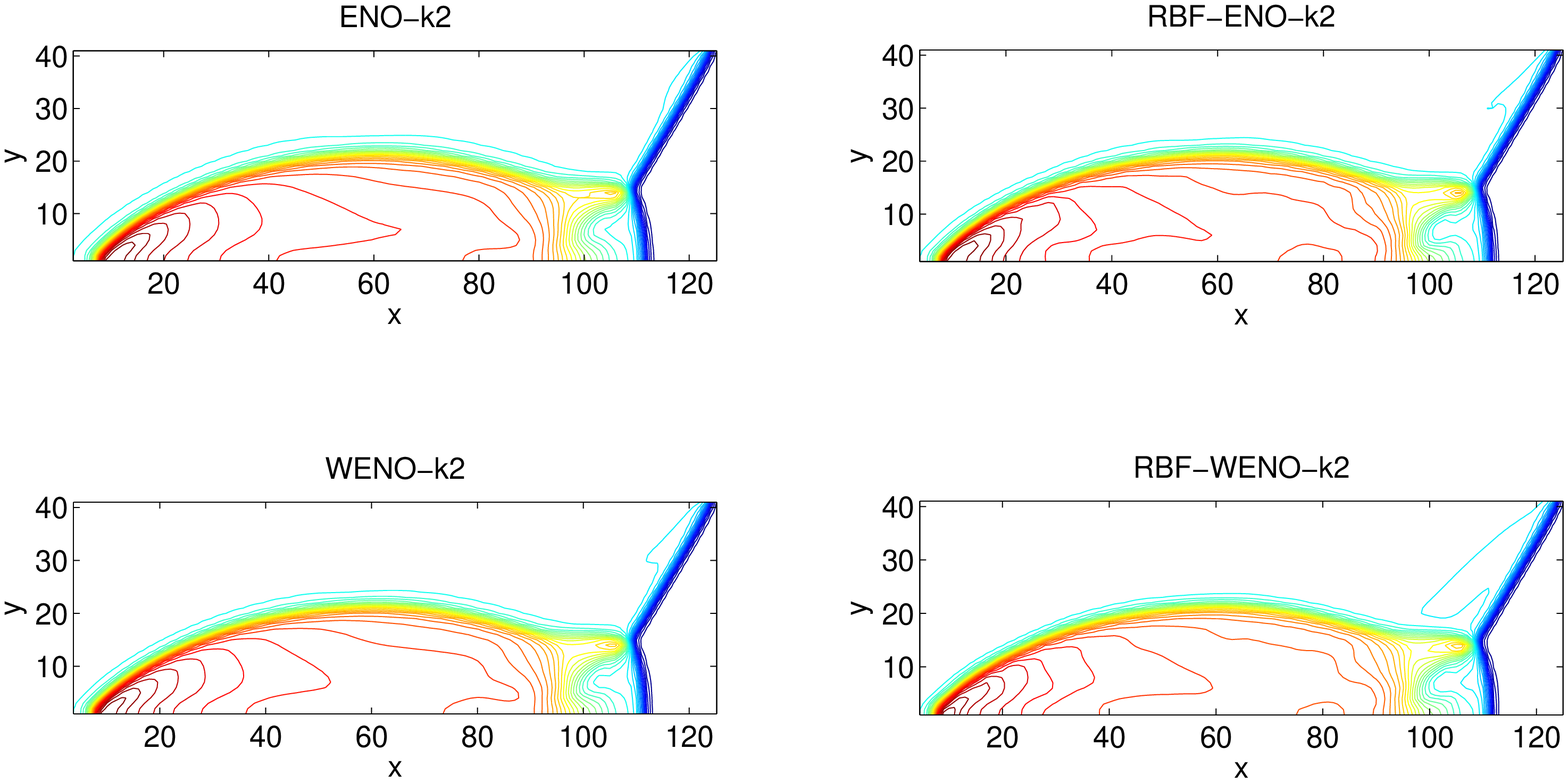}
\hskip -0.6in
\includegraphics[width=0.8\textwidth]{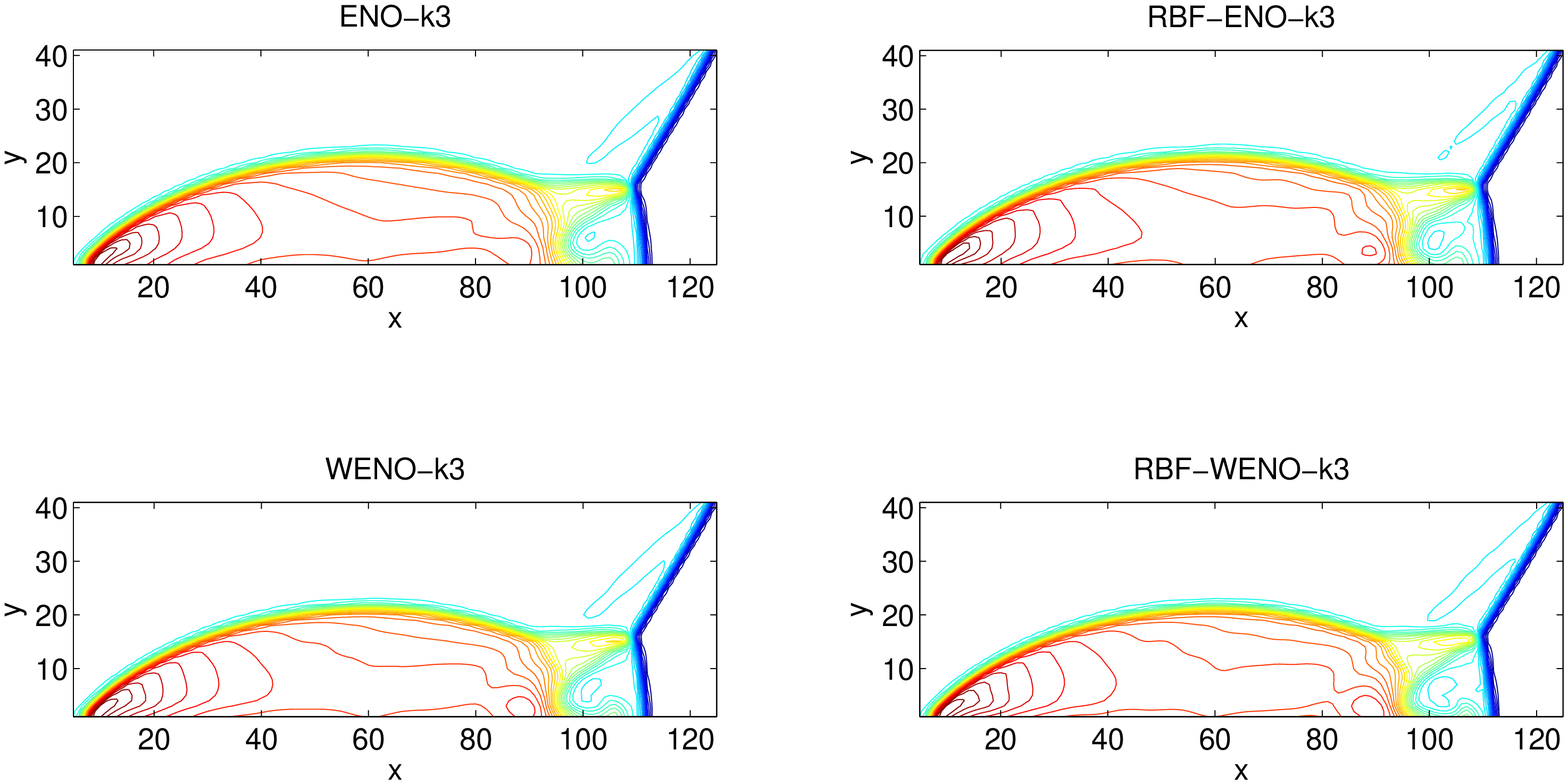}
\end{center}
\caption{Density contours  for the Double Mach Reflection problem. $T=0.2$. $N=160$. $M=40$.}
\label{DMR1}

\bigskip

\begin{center}
\includegraphics[width=0.8\textwidth]{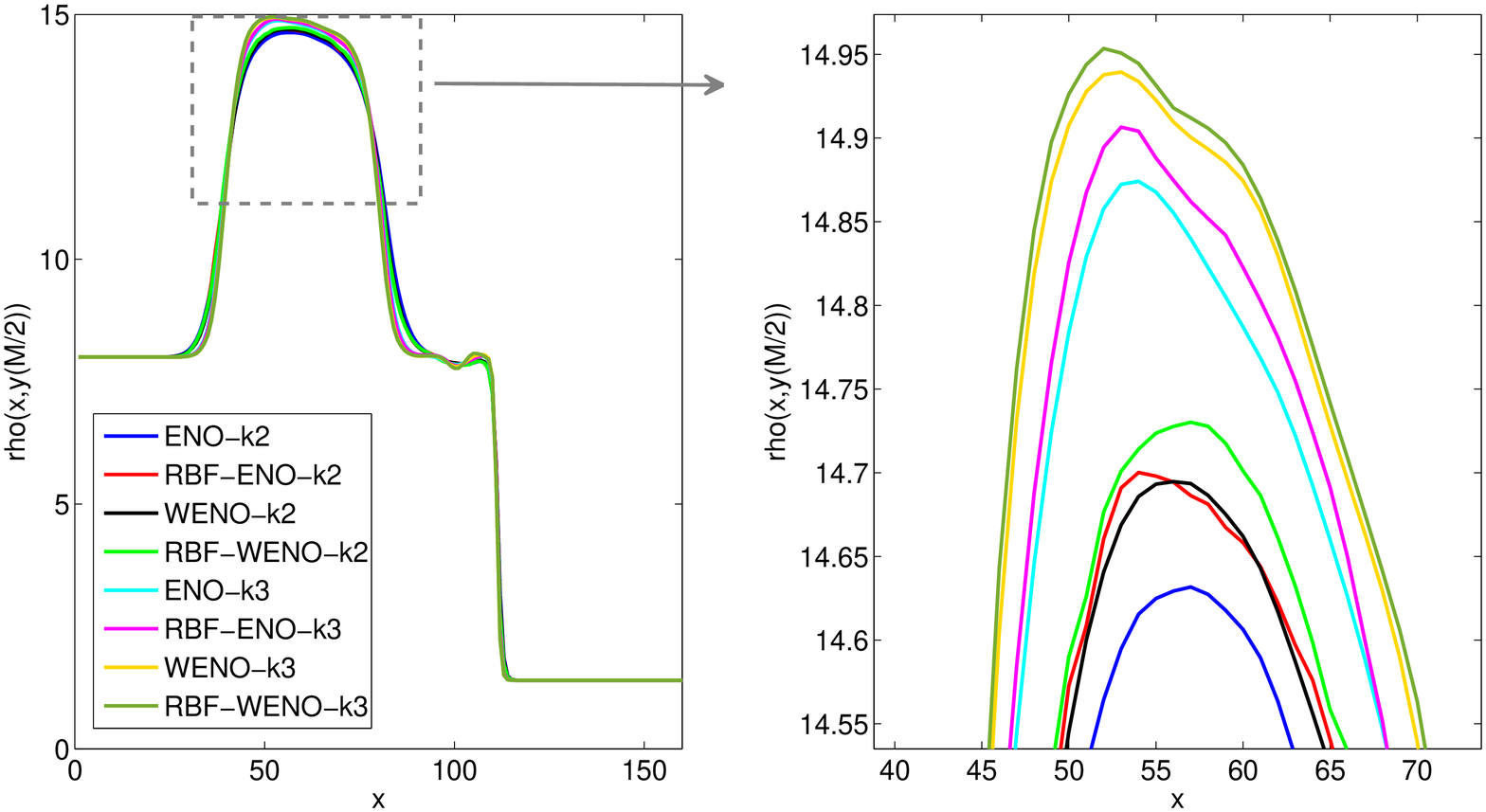}
\end{center}
\caption{Density profiles at $y = 0.5$ for the Double Mach Reflection problem. $T=0.2$. $N=160$. $M=40$. }
\label{DMR2}
\end{figure}
\begin{figure}[h]
\begin{center}
\includegraphics[width=0.8\textwidth]{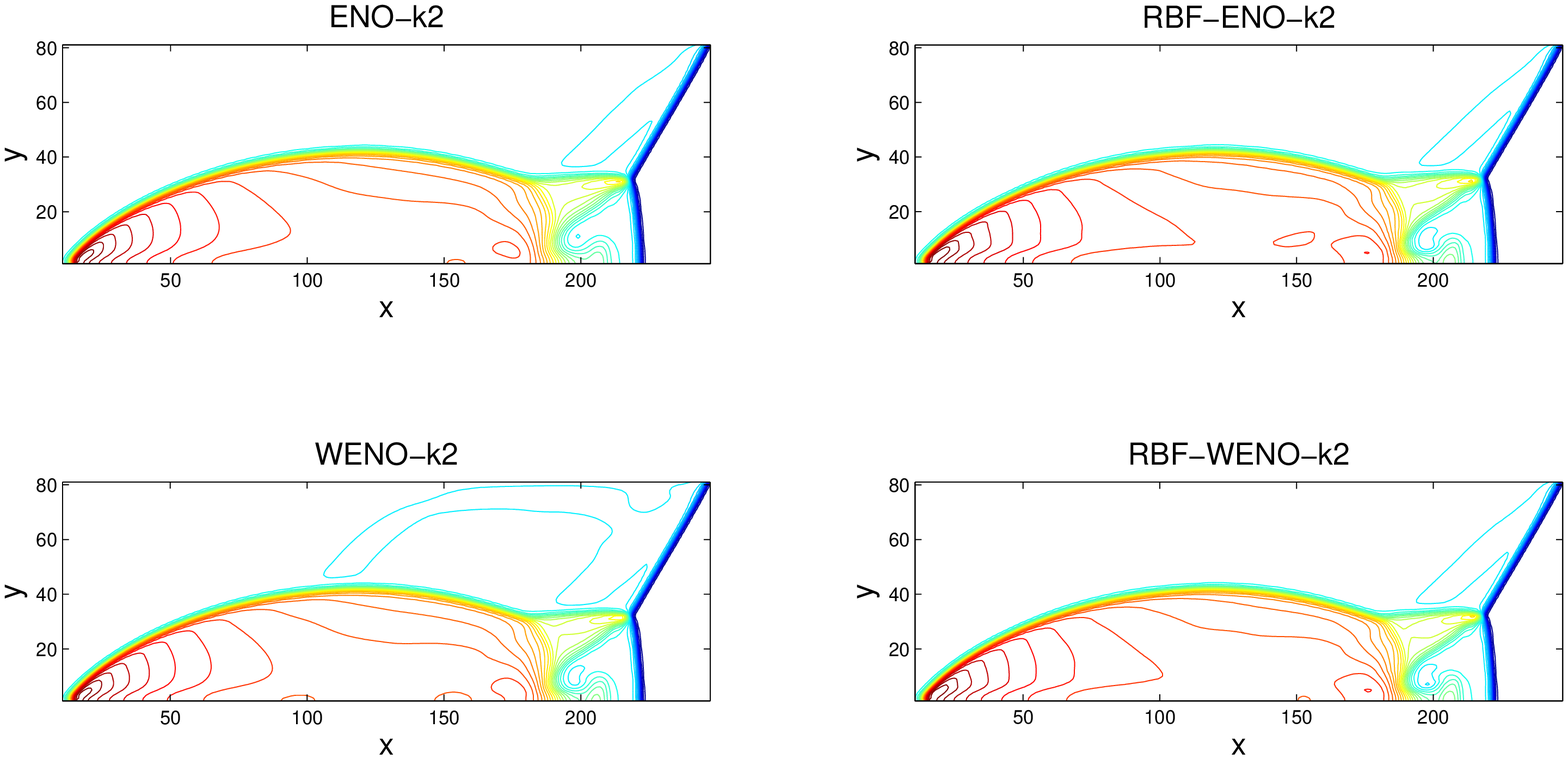}
\hskip -0.6in
\includegraphics[width=0.8\textwidth]{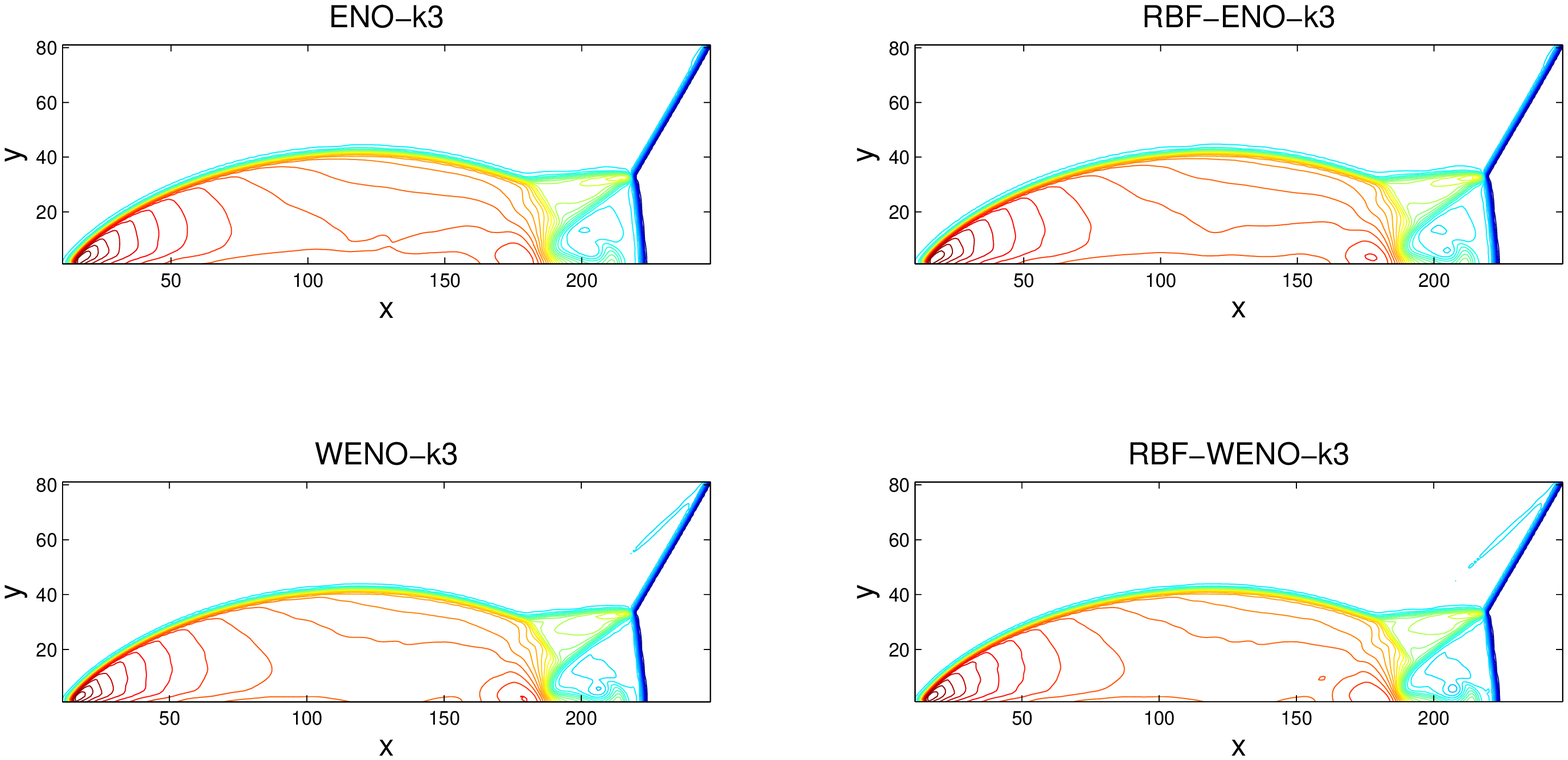}
\end{center}
\caption{Density contours for the Double Mach Reflection problem. $T=0.2$. $N=320$. $M=80$.}
\label{DMR11}
\end{figure}


\begin{figure}[h]
\begin{center}
\includegraphics[width=0.8\textwidth]{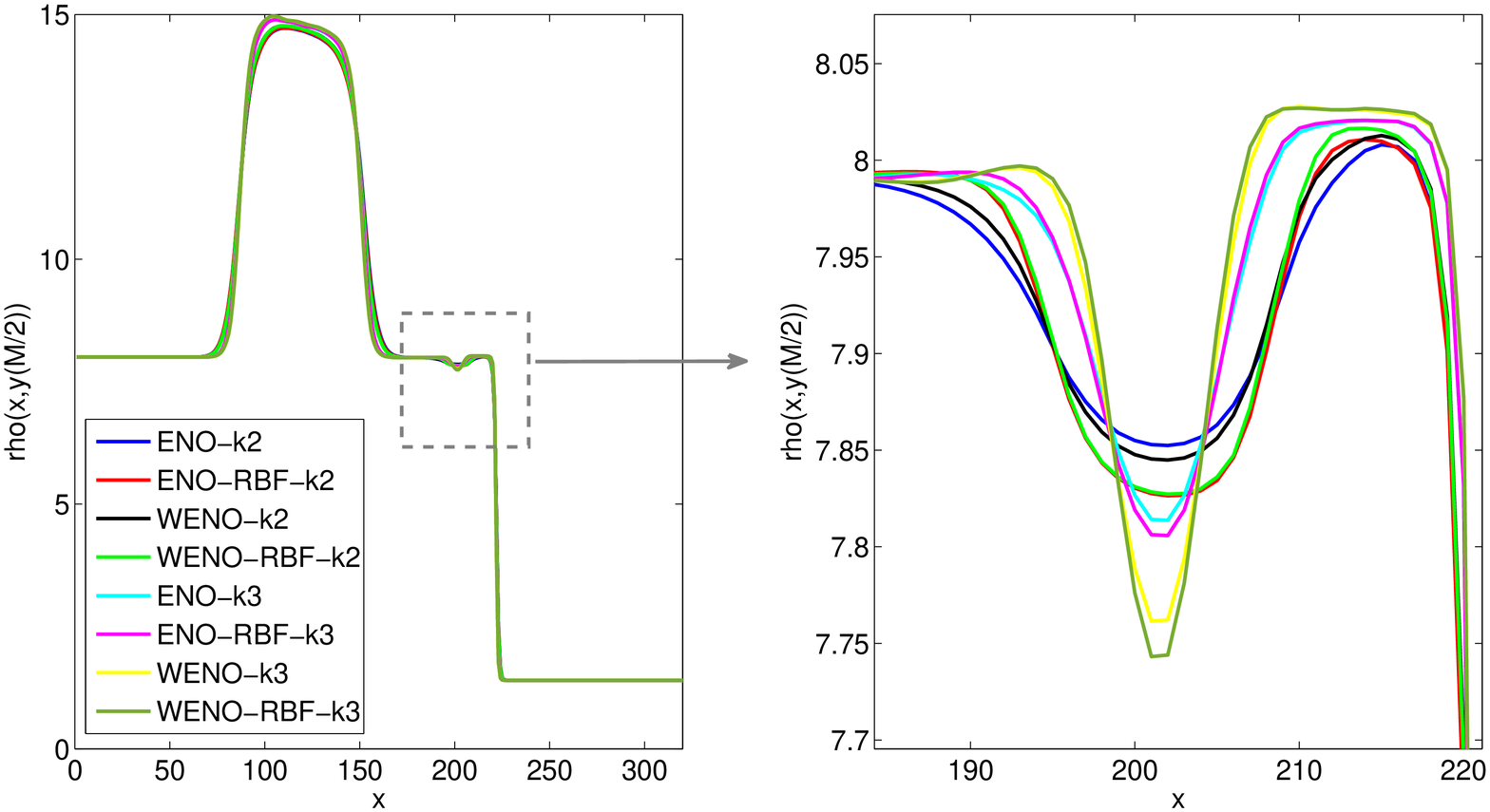}
\end{center}
\caption{Density profiles at $y = 0.5$ for the Double Mach Reflection problem. $T=0.2$. $N=320$. $M=80$. }
\label{DMR22}
\end{figure}

\section{Conclusion}
\label{sec:conclusion}
In this paper, we presented the finite difference RBF-ENO/WENO methods, a direct extension of the non-polynomial ENO/WENO finite volume methods proposed in \cite{GuoJung}. The RBF-ENO/WENO finite difference method seeks the reconstruction using the RBF interpolation which involves the free shape parameter. By optimizing the shape parameter, the finite difference RBF-ENO/WENO reconstruction becomes more accurate than the regular ENO/WENO reconstruction for smooth problems and yields sharper solution profiles near the jump discontinuity. The finite difference RBF-ENO reconstruction is sought in terms of the flux function. Accordingly the optimization of the shape parameter of the RBFs is found based on the flux functions on the grid points within the given stencil. Unlike the finite volume RBF-ENO/WENO method, the finite difference RBF-ENO/WENO method does not involve any integral approximation. Thus the order of convergence is fully determined by the size of stencil for 2D problems. 
Numerical examples both in 1D and 2D show that the finite difference RBF-ENO and RBF-WENO-JS methods yield more accurate results than the regular ENO and WENO-JS methods. Also the solutions by the RBF-ENO and RBF-WENO-JS methods near the discontinuities are sharper than the solutions by the regular ENO/WENO method. 

In our current work, we implemented the RBF reconstruction in the regular ENO and WENO methods based on the original WENO reconstruction by Jiang and Shu \cite{WENO}. As mentioned in the paper, we will investigate the RBF reconstruction with other ENO/WENO variations in our future research such as the WENO-Z method \cite{WENOZ}. We were interested in improving local accuracy, but it would be an interesting research to investigate how the local adaption of the shape parameter can be utilized with the meshless properties of RBFs on the unstructured grid. Furthermore, we will also investigate the RBF-ENO/WENO method for higher values of $k$ than $k = 3$ with multiple shape parameters defined in the RBF basis used for the reconstruction. 

\vskip .1in
\noindent
{\bf Acknowledgments:} The authors thank W.-S. Don for his useful communication. 

\end{document}